\documentclass{amsart}

\usepackage{amsmath} 
\usepackage{amssymb}
\usepackage{mathrsfs}

\newtheorem{theorem}{Theorem}[section] 
\newtheorem{claim}[theorem]{Claim}

\theoremstyle{definition}
\newtheorem{definition}[theorem]{Definition}
\newtheorem{defclaim}[theorem]{Definition/Claim}

\newtheorem{convention}[theorem]{Convention}
\newtheorem{choice}[theorem]{Choice}
\newtheorem{discussion}[theorem]{Discussion}
\newtheorem{problem}[theorem]{Problem}
\newtheorem{observation}[theorem]{Observation} 
\newtheorem{question}[theorem]{Question}

\newtheorem{thesis}[theorem]{Thesis}

\theoremstyle{remark}
\newtheorem{remark}[theorem]{Remark}

\newtheorem{context}[theorem]{Context}

\newcommand{\rest}{{\restriction}}

\newcommand{\wilog}{{\rm without loss of generality}}

\newcommand{\then}{{\underline{then}}}
\newcommand{\when}{{\underline{when}}}
\newcommand{\Then}{{\underline{Then}}}

\newcommand{\mn}{{\medskip\noindent}}
\newcommand{\sn}{{\smallskip\noindent}}

\newcommand{\gb}{{\mathfrak b}}
\newcommand{\gB}{{\mathfrak B}}
\newcommand{\gC}{{\mathfrak C}}

\newcommand{\bbG}{{\mathbb G}}

\newcommand{\cD}{{\mathcal D}}
\newcommand{\gd}{{\mathfrak d\/}} 

\newcommand{\cH}{{\mathcal H}}

\newcommand{\cF}{{\mathcal F}}

\newcommand{\cI}{{\mathcal I}}

\newcommand{\bbP}{{\mathbb P}}
\newcommand{\bbR}{{\mathbb R}}
\newcommand{\cP}{{\mathcal P}}
\newcommand{\gp}{{\mathfrak p}}

\newcommand{\bbQ}{{\mathbb Q}}

\newcommand{\gt}{{\mathfrak t}} 
\newcommand{\cU}{{\mathcal U}}

\newcount\skewfactor
\def\mathunderaccent#1#2 {\let\theaccent#1\skewfactor#2
\mathpalette\putaccentunder}
\def\putaccentunder#1#2{\oalign{$#1#2$\crcr\hidewidth
\vbox to.2ex{\hbox{$#1\skew\skewfactor\theaccent{}$}\vss}\hidewidth}}
\def\name{\mathunderaccent\tilde-3 }

\newenvironment{PROOF}[2][\proofname.]
   {\begin{proof}[#1]\renewcommand{\qedsymbol}{\textsquare$_{\rm #2}$}}
   {\end{proof}}

\begin{document}

\title {Large continuum, oracles \\
 Sh895}
\author {Saharon Shelah}
\address{Einstein Institute of Mathematics\\
Edmond J. Safra Campus, Givat Ram\\
The Hebrew University of Jerusalem\\
Jerusalem, 91904, Israel\\
 and \\
 Department of Mathematics\\
 Hill Center - Busch Campus \\ 
 Rutgers, The State University of New Jersey \\
 110 Frelinghuysen Road \\
 Piscataway, NJ 08854-8019 USA}
\email{shelah@math.huji.ac.il}
\urladdr{http://shelah.logic.at}
\thanks{Research supported by the United States-Israel Binational
 Science Foundation (Grants No. 2002303 and 2006108).
The author thanks Alice Leonhardt for the beautiful typing.
  First Typed - 07/Feb/6} 


\date{December 22, 2009}

\begin{abstract}
Our main theorem is about iterated forcing for making the
continuum larger than $\aleph_2$.  
We present a generalization of \cite{Sh:669} which deal with oracles for
random, (also for other cases and generalities), by 
replacing $\aleph_1,\aleph_2$ by $\lambda,\lambda^+$
(starting with $\lambda = \lambda^{< \lambda} > \aleph_1$).  Well, 
we demand absolute c.c.c.  So we get, e.g. the continuum
is $\lambda^+$ but we can get cov(meagre) $=\lambda$ and we 
give some applications.  
As in \cite{Sh:669}, it is a
``partial" countable support iteration but it is c.c.c.
\end{abstract}

\maketitle
\numberwithin{equation}{section}
\setcounter{section}{-1}

\section{Introduction}

Starting, e.g. with $\bold V \models$ G.C.H. and $\lambda =
\lambda^{< \lambda} > \aleph_1$, we construct a forcing notion 
$\bbP$ of cardinality $\lambda^+$, by a partial of CS iteration but
the result is a c.c.c. forcing.

The general iteration theorems (treated in \S1) seem generally
suitable for constructing universes with MA$_{< \lambda} +
2^{\aleph_0} = \lambda^+$, and taking more care, we should be able to 
get universes without MA$_{< \lambda}$, see \ref{h.21} below.

Our method is to immitate \cite{Sh:669}; concerning the differences, some
are inessential: using games not using diamonds in the framework
itself, (inessential means that we could have in \cite{Sh:669} 
immitate the choice here and vice versa).

An essential difference is that we deal here with large continuum -
$\lambda^+$; we concentrate on the case we shall (in $\bold V^{\Bbb P}$)
have MA$_{< \lambda}$ but e.g. non(null) $= \lambda$ and
 ${\frak b} = \lambda^+$ (or ${\frak b} = \lambda$).

It seems to us that generally:
\begin{thesis}
\label{h.4}  
The iteration theorem here is enough to get results
parallel to known results with $2^{\aleph_0} = \aleph_2$ replacing
$\aleph_1,\aleph_2$ by $\lambda,\lambda^+$.  
\end{thesis}

\noindent
To test this thesis we have asked Bartoszy\'nski to suggest test
problems for this method and he suggests:
\begin{problem}
\label{h.7}  Prove the consistency of each of the
\mn
\begin{enumerate}
\item[$(A)$]   $\aleph_1 < \lambda < 2^{\aleph_0}$ and the
$\lambda$-Borel conjecture, i.e. $A \subseteq {}^\omega 2$ is of
strong measure zero iff $|A| < \lambda$
\sn
\item[$(B)$]   $\aleph_1 <$ non(null) $< 2^{\aleph_0}$, see \ref{bt.7}
\sn
\item[$(C)$]   $\aleph_1 < {\frak b} = \lambda < 2^{\aleph_0}$ the
dual $\lambda$-Borel conjecture (i.e. $A \subseteq {}^\omega 2$ is
strongly meagre iff $|A| < \lambda$)
\sn
\item[$(D)$]   $\aleph_1 < {\frak b} = \lambda < 2^{\aleph_0} +$ the
dual $2^{\aleph_0}$-Borel conjecture
\sn
\item[$(E)$]   combine (A) and (C) and/or combine (A) and (D).
\end{enumerate}
\end{problem}

\noindent
Parallely Steprans suggests:
\begin{problem}
\label{h.14} 
1) Is there a set $A \subseteq {}^\omega 2$ of
cardinality $\aleph_2$ of $p$-Hausdorff measure $>0$, but for every set of
size $\aleph_2$ is null (for the Lebesgue measure)?

\noindent
2) The (basic product) I think ${\gb} = {\gd} \vee {\gd} =
2^{\aleph_0}$ gives an answer, what about cov(meagre) $=\lambda <
2^{\aleph_0}$?
\medskip

We shall deal with the iteration in \S1, give an application to a
problem from \cite{Sh:885} in \S2 (and \S3,\S4).

Lastly, in \S5 we deal with Bartoszy\'nski's test problem (B), in
fact, we get quite general such results.
\end{problem}

\noindent
It is natural to ask
\begin{discussion}
\label{h.21}  
1) In \S1, we may wonder if we can give
``reasonable" sufficient condition for ${\frak b} = \aleph_1$ or
${\frak b} = \kappa < \lambda$?  The answer is yes.  It is natural to assume
that we have in $\bold V$ a 
$<_{J^{\text{bd}}_\omega}$-increasing sequence $\bar f = \langle
f_\alpha:\alpha < \kappa\rangle$ of functions from ${}^\omega \omega$
with no $<^*_{J^{\text{bd}}_\omega}$-upper bound and we would like to
preserve this property of $\bar f$, i.e. in \S1 we
\mn
\begin{enumerate}
\item[$(a)$]  restrict ourselves to $\bold p \in K^1_\lambda$ such 
that $\Vdash_{\Bbb P_{\bold p}} ``\bar f$ as above".
\end{enumerate}
\mn
More formally redefine $K^1_\lambda$ such that
\begin{enumerate}
\item[$(b)$]  replace ``$\bbP$ is absolute c.c.c." by ``$\bbP$ is
  c.c.c., preserve $\bar f$ as above and if $\bbQ$ satisfies those two
  conditions then also $\bbP \times \bbQ$ too satisfies those two
  conditions.
\end{enumerate}
\mn
This has similar closure properties, that is, the proofs do not really
change.

\noindent
2) More generally consider $K$, a property of forcing notions such that:
\mn
\begin{enumerate}
\item[$(a)$]   $\bbP \in K \Rightarrow \bbP$ is c.c.c.
\sn
\item[$(b)$]  $K$ is closed under $\lessdot$-increasing continuous unions
\sn
\item[$(c)$]   $K$ is closed under composition
\sn
\item[$(d)$]   we replace in \S1 ``$\bold p \in K^1_\lambda$" by
``$\bold p \in K$ has cardinality $< \lambda$"
\sn
\item[$(e)$]  we replace in \S1, ``$\bbP$ is absolutely c.c.c." by
``$\bbP \in K$ and $\bbR \in K \Rightarrow \bbP \times \bbR \in K$".
\end{enumerate}
\mn
3) What about using ${\cP}(n)$-amalgamation of forcing notions?  If
   we fix $n$ this seems a natural way to get non-equality for many
   $n$-tuples of cardinal invariants; hopefully we shall return to
this sometime.

\noindent
4) What about forcing by the set of approximations $\bold k$?  See
\ref{it.21}. 
\end{discussion}

\begin{definition}
\label{h.35}  
1) We say a forcing notion $\bbP$ is absolutely c.c.c. \when  \, 
for every c.c.c. forcing notion $\bbQ$ we have $\Vdash_{\bbQ} ``\bbP$ 
is c.c.c."

\noindent
2) We say $\bbP_2$ is absolutely c.c.c. over $\bbP_1$ \when \,
($\bbP_1 \lessdot \bbP_2$ and) $\bbP_2/\bbP_1$ is absolutely c.c.c.  

\noindent
3) Let $\bbP_1 \subseteq_{\text{ic}} \bbP_2$ means that $\bbP_1
\subseteq \bbP_2$ (as quasi orders) and if $p,q \in \bbP_1$ are
incompatible in $\bbP_1$ then they are incompatible in $\bbP_2$
(the inverse holds too) .

The following tries to describe the iteration theorem, this may be
more useful to the reader after having a first reading of \S1.

We treat $\lambda$ as the vertical direction and
$\lambda^+$ as the horizontal direction, the meaning will be clarified
in \S2; our forcing is the increasing union of 
$\langle \bbP^{\bold k_\varepsilon}:\varepsilon < 
\lambda^+\rangle$ where $\bold k_\varepsilon \in K_2$ (so $\bold
k_\varepsilon$ gives an iteration $\langle \bbP_\alpha[\bold
k_\varepsilon]:\alpha <  \lambda\rangle$, 
i.e. a $\lessdot$-increasing continuous sequence of c.c.c. forcing notions)
 and for each such $\bold k_\varepsilon$ each iterand $\Bbb P_{\bold
 p_\alpha[\bold k_\varepsilon]}$ is of cardinality $< \lambda$ and for
 each $\varepsilon < \lambda^+$ the forcing notion
$\Bbb P^{\bold k_\varepsilon}$ is the union of increasing union of
 continuous sequence $\langle \Bbb P_{\bold p_\alpha[\bold k_\varepsilon]}:
\alpha < \lambda\rangle$.  So we can say that $\Bbb P^{\bold
k_\varepsilon}$ is the limit of an FS iteration of length 
$\lambda$, each iterand of cardinality $< \lambda$ and for $\zeta \in
(\varepsilon,\lambda^+),\bold k_\zeta$ gives a ``fatter" iteration,
which for ``most" $\delta \in S (\subseteq \lambda)$, is a 
reasonable extension.  
\end{definition}

\begin{question}
\label{h.42}
Can we get something interesting for the 
continuum $> \lambda^+$ and/or get cov(meagre) $< \lambda$?  This
certainly involves some losses!  We intend to try elsewhere.
\end{question}

\begin{definition}
\label{h45}  
1) For a set $x$ let otrcl$(x)$, the transitive closure over the
   ordinals of $x$, be the minimal set $y$ such that $x \in y 
\wedge (\forall t \in y)(t \notin \text{ Ord } \rightarrow t \subseteq y]$.

\noindent
2) For a set $u$ of ordinals let ${\cH}_{< \kappa}(u)$ be the set
   of $x$ such that otrcl$(x)$ is a subset of $u$ or cardinality $< \kappa$.
\end{definition}

\begin{remark}
\label{h48}  
0) We use $\cH_{< \kappa}(u)$ (in Definition \ref{it.7}) just 
for bookkeeping convenience.

\noindent
1) It is natural to have Ord, the class of ordinals, a class of
   urelements.

\noindent
2) If $\omega_1 \subseteq u$ for ${\cH}_{< \aleph_1}(u)$ it makes
   no difference, but if $\omega_1 \nsubseteq u$ and $\beta = \text{
   min}(\omega_1 \backslash u)$ \then \, $\beta$ is a countable
 subset of $u$ but $\notin {\cH}_{< \aleph_1}(u)$.  Also we use
   ${\cH}_{< \aleph_0}(u)$ where $\omega \subseteq u$, so there are no
   problems. 
\end{remark}
\newpage

\section {The iteration theorem}

If we use the construction for $\lambda = \aleph_1$, the version we
get is closer to, but not the same as \cite{Sh:669} with the 
forcing being locally Cohen.

Here there are ``atomic" forcings used below coming from three sources: 
\mn
\begin{enumerate}
\item[$(a)$]   the forcing given by the winning strategies $\bold
s_\delta$(see below), i.e. the quotient
\sn
\item[$(b)$]   forcing notions intended to generate MA$_{< \lambda}$
\newline
[see \ref{it.56}; we are given $\bold k_1 \in K^2_f$, an approximation
of size $\lambda$, see Definition \ref{it.21}, and a 
$\bbP_{\bold k_1}$-name $\name{\bbQ}$ of a
c.c.c. forcing and sequence $\langle \name{\cI}_i:
i < i(*)\rangle$ of $< \lambda$ dense subsets of $\bbQ$.  We
would like to find $\bold k_2 \in K_2$ satisfying 
$\bold k_1 \le_{K^2_f} \bold k_2$
such that $\Vdash_{\bbP_{\bold k_2}}$ ``there is a directed $G
\subseteq \name{\bbQ}$ not disjoint to any
$\name{\cI}_i (i < i(*))$".  We \underline{do not}
use composition, only $\bbP_{\bold p_\alpha[\bold k_2]} = 
\bbP_{\bold p_\alpha[\bold k_1]} * 
\name{\bbQ}$ for some $\alpha \in E_{\bold k_1} \cap E_{\bold k_2}]$]
\sn
\item[$(c)$]   given $\bold k_1 \in K^2_f$, and
$\name{\bbQ}$ which is a $\bbP_{\bold k_1}$-name of a suitable c.c.c. 
forcing of cardinality $\lambda$ can we find $\bold k_2$ such that
$\bold k_1 \le_{K^2_f} \bold k_2$ and in $\bold V$ we have 
$\Vdash_{\bbP[\bold k_2]}$ 
``there is a subset of $\name{\bbQ}$ generic over $\bold V[\name G
\cap \bbP_{\bold k_1}]"$.
\end{enumerate}
\mn
Let us describe the roles of some of the definitions.  We shall
construct (in the main case) a forcing notion of cardinality
$\lambda^+$ by approximations $\bold k \in K^2_f$ of size (= cardinality)
$\lambda$, see Definition \ref{it.21}, which 
are constructed by approximations $\bold p \in K_1$ of
cardinality $< \lambda$, see Definition \ref{it.7}.

Now $\bold p \in K_1$ is essentially a forcing notion of cardinality 
$< \lambda$, i.e. $\bbP_{\bold p} = 
(P_{\bold p},\le_{\bold p})$, and we add the set $u = u_{\bold p}$
to help the bookkeeping, so (in the main case) $u_{\bold p} 
\in [\lambda^+]^{< \lambda}$.  For the 
bookkeeping we let $P_{\bold p} \subseteq {\cH}_{< \aleph_1}
(u_{\bold p})$, see \ref{h45}(2).

More specifically $\bold k$ (from Definition \ref{it.21}) is mainly a
$\lessdot$-increasing continuous sequence $\bar{\bold p} = \langle
\bold p_\alpha:\alpha \in E_{\bold k}\rangle = \langle \bold
p_\alpha[\bold k]:\alpha \in E_{\bold k}\rangle$, where $E_{\bold k}$ is
a club of $\lambda$.  Hence $\bold k$ represents the forcing notion 
$\bbP_{\bold k} = \cup\{(P_{\bold p_\alpha},\le_{\bold p_\alpha}):\alpha <
\lambda\}$; the union of a $\lessdot$-increasing continuous sequence
of forcing notions $\bbP_{\bold p_\alpha} = \bbP[\bold p_\alpha] =
(P_{\bold p_\alpha},\le_{\bold p_\alpha})$, so we
can look at $\bbP_{\bold k}$ as a FS-iteration.  But then we would
like to construct say an ``immediate successor" $\bold k^+$ of $\bold
k$, so in particular  $\bbP_{\bold k} \lessdot \bbP_{\bold k^+}$, 
e.g. taking care of (b) above so $\name{\bbQ}$ is a 
$\bbP_{\bold k}$-name and even a $\bbP_{\text{min}(E_{\bold k})}$-name
of a c.c.c. forcing notion.
Toward this we choose $\bold p^{\bold k^+}_\alpha = p_\alpha[\bold
k^+]$ by induction on
$\alpha \in E_{\bold k}$.  So it makes sense to demand $\bold p_\alpha
\le_{K_1} \bold p_\alpha[\bold k^+]$, which naturally implies that
$u[\bold p_\alpha] \subseteq u[\bold p^{\bold k^+}_\alpha],
\bbP_{\bold p_\alpha} \lessdot \bbP_{\bold p_\alpha[\bold k^+]}$.
So as $\bold p_\alpha[\bold k^+]$ for $\alpha \in E_{\bold k}$ is
$\le_{K_1}$-increasing continuous, the main case is when
$\beta = \text{ min}(E_{\bold k} \backslash (\alpha +1)$, can we
choose $\bold p_\beta[\bold k^+]$?  
\bigskip

\noindent
Let us try to draw the picture:
\smallskip

\centerline {$\bbP_{\bold p_\beta[\bold k]} \qquad \dashrightarrow \qquad$  ?} 
\smallskip

\centerline {$\uparrow  \hskip75pt \uparrow$}
\smallskip

\centerline {$\bbP_{\bold p_\alpha[\bold k]} \qquad 
\overset\lessdot {}\to \longrightarrow \qquad 
\Bbb P_{\bold p_\alpha[\bold k^+]}$} 
\mn
So we have three forcing notions,
$\bbP_{\bold p_\alpha[\bold k]},\bbP_{\bold p_\beta[\bold k]},
\bbP_{\bold p_\alpha[\bold k^+]}$, where the second and third are
$\lessdot$-extensions of the first.  The main problem is the c.c.c.  As
in the main case we like to have MA$_{< \lambda}$,  there is no
restriction on $\bbP_{\bold p_\alpha[\bold k^+]}/
\bbP_{\bold p_\alpha[\bold k]}$, so it is natural to demand ``$\bbP_{\bold
p_\beta[\bold k]}/\bbP_{\bold p_\alpha[\bold k]}$ is absolutely
c.c.c. for $\alpha < \beta$ from $E_{\bold k}$" 
(recall $\bold p_\alpha[\bold k]$ is demanded to be $<^+_{K_1}$-increasing).

How do we amalgamate?  There are two natural ways which say that ``we
leave $\bbP_{\bold p_\beta[\bold k]}/\bbP_{\bold p_\alpha[\bold k]}$ as it is".
\medskip

\noindent
\underline{First way}:  We decide that $\bbP_{\bold p_\beta[\bold k^+]}$ is
$\bbP_{\bold p_\alpha[\bold k]} \ast ((\bbP_{\bold p_\alpha[\bold k^+]}/
\bbP_{\bold p_\alpha[\bold k]}) \times (\bbP_{\bold p_\beta[\bold k]}
/\bbP_{\bold p_\alpha[\bold k]}))$.

[This is the ``do nothing" case, the lazy man strategy, which  
in glorified fashion we may say: do nothing when in doubt.  Note that
$\bbP_{\bold p_\alpha[\bold k^+]} /
\bbP_{\bold p_\alpha[\bold k]}$ and $\bbP_{\bold p_\beta[\bold k]}/
\bbP_{\bold p_\alpha[\bold k]}$ are $\bbP_{\bold p_\alpha[\bold
k]}$-names of forcing notions.]
\medskip

\noindent
\underline{Second way}:  $\bbP_{\bold p_\beta[\bold k]}/
\bbP_{\bold p_\alpha[\bold k]}$ is defined in some way, e.g. is
a random real forcing in the universe 
$\bold V[\bbP_{\bold p_\alpha[\bold k]}]$ and we decide 
that $\bbP_{\bold p_\beta[\bold k^+]}/\bbP_{\bold p_\alpha[\bold
k^+]}$ is defined in the same way:
the random real forcing in the universe 
$\bold V[\bbP_{\bold p_\alpha[\bold k^+]}]$; this is expressed by
the strategy $\bold s_\alpha$.

[That is: retain the same definition of the forcing in the $\alpha$-th
place, so in some sense we again do nothing novel.]

\begin{context}
\label{it.1}  
Let $\lambda = \text{ cf}(\lambda) > \aleph_1$ or just\footnote{if
$\lambda = \aleph_1$, we can change the definitions of $\bold k \in
K_2$, instead $\langle \bbP_\alpha[\bold k]:\alpha < \lambda\rangle$
is $\lessdot$-increasing, we carry with us large enough family of
dense subsets, e.g. coming from some countable $N$.} $\lambda = \text{
cf}(\lambda) \ge \aleph_1$.
\end{context}

Below, $\le^+_{K_1}$ is used in defining $\bold k \in K^2_f$ as
consisting also of $\le^+_{K_1}$-increasing continuous sequence $\langle
\bold p_\alpha:\alpha \in E \subseteq \lambda\rangle$ (so increasing
vertically). 
\begin{definition}
\label{it.7}  
1) Let $K_1$ be the class of $\bold p$ such that:
\mn
\begin{enumerate}
\item[$(a)$]   $\bold p = (u,P,\le) = (u_{\bold p},P_{\bold p},
\le_{\bold p}) = (u_{\bold p},\Bbb P_{\bold p})$
\sn
\item[$(b)$]   $\omega \subseteq u \subseteq \text{ Ord}$,
\sn
\item[$(c)$]   $P$ is a set $\subseteq {\cH}_{< \aleph_1}(u)$,
\sn
\item[$(d)$]   $\le$ is a quasi-order on $P$,
\end{enumerate}
\mn
satisfying
\mn
\begin{enumerate}
\item[$(e)$]   the pair $(P,\le)$ which we denote also by $\bbP =
\bbP_{\bold p}$ is a c.c.c. forcing notion.
\end{enumerate}
\mn
1A) We may write $u[\bold p],P[\bold p],\Bbb P[\bold p]$.

\noindent
2) $\le_{K_1}$ is the following two-place relation on $K_1:\bold p \le_{K_1}
\bold q$ iff $u_{\bold p} \subseteq u_{\bold q}$ and $\bbP_{\bold p} 
\lessdot \bbP_{\bold q}$ and $\bbP_{\bold q} \cap {\cH}_{< \aleph_1}
(u_{\bold p}) = \bbP_{\bold p}$; moreover, just for transparency $q
\le_{\bbP[\bold q]} p \in \bbP_{\bold p} \Rightarrow q \in \bbP_{\bold p}$.  

\noindent
3) $\le^+_{K_1}$ is the following 
two-place relation on $K_1:\bold p \le^+_{K_1}
\bold q$ \underline{iff} $\bold p \le_{K_1} \bold q$ and 
$\bbP_{\bold q}/\bbP_{\bold p}$ is absolutely c.c.c., see
Definition \ref{h.35}(1).

\noindent
4) $K^1_\lambda$ is the family of $\bold p \in K_1$ such that
 $u_{\bold p} \subseteq \lambda^+$ and $|u_{\bold p}| < \lambda$.

\noindent
5) We say $\bold p$ is \underline{the exact limit} 
of $\langle \bold p_\alpha:\alpha \in v\rangle,v 
\subseteq \text{ Ord}$, in symbols $\bold p =
\cup\{\bold p_\alpha:\alpha \in v\}$ when $u_{\bold p} =
\cup\{u_{\bold p_\alpha}:\alpha \in v\},\bbP_{\bold p} =
 \cup\{\bbP_{\bold p_\alpha}:\alpha \in v\}$ and $\alpha \in v 
\Rightarrow \bold p_\alpha \le_{K_1} \bold p$; hence $\bold p \in K_1$.

\noindent
6) We say $\bold p$ is just \underline{a limit} of $\langle 
\bold p_\alpha:\alpha \in v\rangle$ when 
$u_{\bold p}$ is $\cup\{u_{\bold p_\alpha}:\alpha \in v\},
\bbP_{\bold p} \supseteq \cup\{\bbP_{\bold p_\alpha}:\alpha \in v\}$ 
and $\alpha \in v \Rightarrow \bold p_\alpha \le_{K_1} \bold p$.

\noindent
7) We say $\bar{\bold p} = \langle \bold p_\alpha:\alpha <
\alpha^*\rangle$ is $\le_{K_1}$-increasing continuous [strictly
$\le_{K_1}$-increasing continuous] \underline{when} it is
$\le_{K_1}$-increasing and for every limit $\alpha < \alpha^*,\bold
p_\alpha$ is a limit of $\bar{\bold p} \restriction \alpha$ [is the exact
limit of $\bar{\bold p} \restriction \alpha$], respectively.
\end{definition}

\begin{observation}
\label{it.3}  
1) $\le_{K_1}$ is a partial order on $K_1$.

\noindent
2) $\le^+_{K_1} \subseteq \le_{K_1}$ is a partial order on $K_1$.

\noindent
3) If $\bar{\bold p} = \langle \bold p_\alpha:\alpha < \delta\rangle$
is a $\le_{K_1}$-increasing sequence and 
$\cup\{\bbP_{\bold p_\alpha}:\alpha < \delta\}$ satisfies the c.c.c.
and $\delta < \lambda$ \then \, some $\bold p \in K_1$ 
is the union $\cup\{\bold p_\alpha:\alpha
< \delta\}$ of $\bar{\bold p}$, i.e. $\cup \bar{\bold p} \in K_1$ and
$\alpha < \delta \Rightarrow \bold p_\alpha \le_{K_1} \bold p$; this
determines $\bold p$ uniquely and $\bold p$ is the exact union of
$\bar{\bold p}$.

\noindent
4) If $\bar{\bold p} = \langle \bold p_\alpha:\alpha < \delta\rangle$
is $\le_{K_1}$-increasing and {\rm cf}$(\delta) = \aleph_1$ implies
$\{\alpha < \delta:\bold p_\alpha$ the exact limit of $\bar{\bold p}
\restriction \alpha$ or just $\bigcup\limits_{\beta < \alpha} 
\bbP_{\bold p_\beta}
\lessdot \bbP_{\bold p_\alpha}\}$ is a stationary subset 
of $\delta$ \then \, $\cup\bar{\bold p} \in K_1$ 
is a $\le_{K_1}$-upper bound of $\bar{\bold p}$ and is the exact limit
of $\bar{\bold p}$.
 
\noindent
5) If in part (4), 
$\bar{\bold p}$ is also $\le^+_{K_1}$-increasing \then \,
 $\alpha < \delta \Rightarrow \bold p_\alpha \le^+_{K_1} \bold p$.
\end{observation}

\begin{PROOF}{\ref{it.3}}
Should be clear, e.g. in part (5) recall that
c.c.c. forcing preserve stationarity of subsets of $\delta$.  
\end{PROOF}

\noindent
We now define the partial order $\le^*_{K_1}$; it will be used in 
describing $\bold k_1 <_{K_2} \bold k_2$, i.e. demanding $(\bold p^{\bold
k_1}_\alpha,\bold p^{\bold k_2}_\alpha) \le^*_{K_1} (\bold p^{\bold
k_1}_{\alpha +1},\bold p^{\bold k_2}_{\alpha +1})$ for many $\alpha <
\lambda$.  
\begin{definition}
\label{it.4} 
1) Let $\le^*_{K_1}$ be the following
two-place relation on the family of pairs 
$\{(\bold p,\bold q):\bold p \le_{K_1} \bold q\}$.  We let
$(\bold p_1,\bold q_1) \le^*_{K_1} (\bold p_2,\bold q_2)$ \underline{iff}
\mn
\begin{enumerate}
\item[$(a)$]   $\bold p_1 \le^+_{K_1} \bold p_2$
\sn
\item[$(b)$]   $\bold q_1 \le^+_{K_1} \bold q_2$
\sn
\item[$(c)$]   $\Vdash_{\bbP[\bold p_2]} ``\bbP_{\bold q_1}/
(\name G_{\bbP[\bold p_2]} \cap \bbP_{\bold p_1}) \lessdot \bbP_{\bold q_2}/
\name G_{\bbP[\bold p_2]}"$
\sn
\item[$(d)$]   $u_{\bold p_2} \cap u_{\bold q_1} = u_{\bold p_1}$ 
\end{enumerate}
\mn
2) Let $\le'_{K_1}$ be the following two-place relation on the family
 $\{(\bold p,\bold q):\bold p \le_{K_1} \bold q\}$ of pairs.  We let
 $(\bold p_1,\bold q_1) \le'_{K_1} (\bold p_2,\bold q_2)$ 
\underline{iff} clauses (a),(b),(d) from part (1) above and
\mn
\begin{enumerate}
\item[$(c)'$]   if $p_1 \in \bbP_{\bold p_1},q_1 \in \bbP_{\bold q_1}$ and
$p_1 \Vdash_{\bbP_{\bold p_1}} ``q_1 \in \bbP_{\bold p_2}/
\name G_{\bbP_{\bold p_1}}"$ \then \, $p_1 \Vdash_{\bbP_{\bold p_2}} 
``q_1 \in \Bbb P_{\bold q_2}/\name G_{\bbP_{\bold p_2}}"$.
\end{enumerate}
\mn
3) Assume $\bold p_\ell \in K_1$ for $\ell=0,1,2$ and $\bold p_0
\le_{K_1} \bold p_1$ and $\bold p_0 \le_{K_1} \bold p_2$ and
$u_{\bold p_1} \cap u_{\bold p_2} = u_{\bold p_0}$.  We define the
 amalgamation $\bold p = \bold p_3 = \bold p_1 \times_{\bold p_0} \bold p_2$ or
 $\bold p_3 = \bold p_1 \times \bold p_2/\bold p_0$ as the triple
 $(u_{\bold p},P_{\bold p} \le_{\bold p})$ as follows\footnote{If in
clause (b) of \ref{it.4}(3) we would like to avoid ``$p_\ell \in
\bbP_{\bold p_\ell} \backslash \bbP_{\bold p_0}"$ we may replace
$(p_1,p_2)$ by $(p_1,p_2,u_{\bold p_1} \cup u_{\bold p_2})$ when
$\bold p_1 \ne \bold p_1 \wedge \bold p_0 \ne \bold p_2$ equivalently
$\bold p_0 \ne \bold p_1 \wedge \bold p_0 \ne \bold p_2$.}:
\mn
\begin{enumerate}
\item[$(a)$]   $u_{\bold p} = u_{\bold p_1} \cup u_{\bold p_2}$
\sn
\item[$(b)$]   $P_{\bold p} = P_{\bold p_1} \cup P_{\bold p_2}
\cup\{(p_1,p_2):p_1 \in P_{\bold p_1} \backslash P_{\bold p_0},
p_2 \in P_{\bold p_2} \backslash P_{\bold p_0}$ and
   for some $p \in P_{\bold p_0}$ we have $p \Vdash_{P[\bold p_0]}
``p_\ell \in \bbP_{p_\ell}/\bbP_{\bold p_0}"$ for $\ell=1,2\}$
\sn
\item[$(c)$]  $\le_{\bold p}$ is defined naturally as $\le_{\bold p_1} 
\cup \le_{\bold p_2} \cup\{((p_1,p_2),(q_1,q_2)):
(p_1,p_2),(q_1,q_2) \in P_{\bold p}$
and $p_1 \le_{\bold p_1} q_1$ and $p_2 \le_{\bold p_2} q_2\} \cup
\{(p'_\ell,(p_1,p_2)):p'_\ell \in P_{\bold p_\ell},(p_1,p_2) \in
P_{\bold p}$ and $p'_\ell \le_{\bold p_1} p_\ell$ and $\ell \in \{1,2\}\}$.
\end{enumerate}
\end{definition}

\begin{remark}
Why not use $u$ instead ${\cH}_{< \aleph_1}(u)$?  Not a real
 difference but, e.g. there may not be enough elements in a union of two.
\end{remark}

\begin{observation}
\label{it.4.7} 
1) $\le^*_{K_1},\le'_{K_1}$ are partial orders on their domains. 

\noindent
2) $(\bold p_1,\bold q_1) \le^*_{K_1} (\bold p_1,\bold q_1)$ 
implies $(\bold p_1,\bold q_1) \le'_{K_1} (\bold p_2,\bold q_2)$.
\end{observation}

\noindent
For the ``successor case vertically and horizontally" we shall use
\begin{claim}
\label{it.8}  
Assume that $\bold p_1 \le^+_{K_1} \bold p_2$ and 
$\bold p_1 \le_{K_1} \bold q_1$ and $u_{\bold p_2} \cap
 u_{\bold q_1} = u_{\bold p_1}$ \underline{then} $\bold q_2 \in K_1$
and $(\bold p_1,\bold q_1)
\le^*_{K_1} (\bold p_2,\bold q_2)$ \when \, we define $\bold q_2 =
\bold q_1 \times_{\bold p_1} \bold p_2$ as in \ref{it.4}(3).
\end{claim}

\begin{PROOF}{\ref{it.8}}
 Straight. 
\end{PROOF}

\noindent
The following claim will be applied to a pair of vertically
increasing continuous sequences, one laying horizontally above the other.
\begin{claim}
\label{it.5.21}  
Assume $\varepsilon(*) < \lambda$ and
\mn
\begin{enumerate}
\item[$(a)$]   $\langle \bold p^\ell_\varepsilon:\varepsilon \le
\varepsilon(*)\rangle$ is $\le^+_{K_1}$-increasing continuous for $\ell=1,2$
\sn
\item[$(b)$]   $(\bold p^1_\varepsilon,\bold p^2_\varepsilon)
\le'_{K_1} (\bold p^1_{\varepsilon +1},\bold p^2_{\varepsilon +1})$ 
for $\varepsilon < \varepsilon(*)$.
\end{enumerate}
\mn
\Then 
\mn
\begin{enumerate}
\item[$(\alpha)$]   $\bold p^1_{\varepsilon(*)} \le_{K_1} \bold
p^2_{\varepsilon(*)}$
\sn
\item[$(\beta)$]   for $\varepsilon < \zeta \le \varepsilon(*)$ we
have $(\bold p^1_\varepsilon,\bold p^2_\varepsilon) \le'_{K_1}
(\bold p^1_\zeta,\bold p^2_\zeta)$.
\end{enumerate}
\end{claim}

\begin{PROOF}{\ref{it.5.21}}
Easy. 
\end{PROOF}

\noindent
For the ``successor case horizontally, limit case vertically when the
relevant game, i.e. the relevant winning strategy is not active" we shall use
\begin{claim}
\label{it.5}  
Assume $\varepsilon(*) < \lambda$ is a limit ordinal and
\mn
\begin{enumerate}
\item[$(a)$]   $\langle \bold p_\varepsilon:\varepsilon \le
\varepsilon(*)\rangle$ is $\le^+_{K_1}$-increasing, and $\langle \bold
q_\varepsilon:\varepsilon < \varepsilon(*)\rangle$ is $\le^+_{K_1}$-increasing
\sn
\item[$(b)$]   $\bold p_\varepsilon \le_{K_1} \bold q_\varepsilon$ 
for $\varepsilon < \varepsilon(*)$
\sn
\item[$(c)$]   if $\varepsilon < \zeta < \varepsilon(*)$ then
$(\bold p_\varepsilon,\bold q_\varepsilon) \le'_{K_1} (\bold
p_\zeta,\bold q_\zeta)$
\sn
\item[$(d)$]   if $\zeta < \varepsilon(*)$ is a limit ordinal then
$\Vdash_{\bbP_{[\bold p_\zeta]}} ``\bbP_{\bold q_\zeta}/
\name G_{\bbP[\bold p_\zeta]} = \cup \{\bbP_{\bold q_\varepsilon}/
(\name G_{\bbP_{[\bold p_\zeta]}} \cap \bbP_{\bold p_\varepsilon}):
\varepsilon < \zeta\}$.
\end{enumerate}
\mn
\Then \, we can choose $\bold q_{\varepsilon(*)}$ such that
\mn
\begin{enumerate}
\item[$(\alpha)$]   $\bold p_{\varepsilon(*)} \le_{K_1} 
\bold q_{\varepsilon(*)}$
\sn
\item[$(\beta)$]   $(\bold p_\varepsilon,\bold q_\varepsilon)
\le'_{K_1} (\bold p_{\varepsilon(*)},\bold q_{\varepsilon(*)})$ for
every $\varepsilon < \varepsilon(*)$
\sn
\item[$(\gamma)$]   clause (d) holds also for $\zeta =
\varepsilon(*)$.
\end{enumerate}
\end{claim}

\begin{remark}
We can replace $\le'_{K_1}$ by $\le^*_{K_1}$ in $(c)$ and $(\beta)$.
\end{remark}

\begin{PROOF}{\ref{it.5}}  
Should be clear.  
\end{PROOF}

\noindent
The game defined below is the non-FS ingredient; 
(in the main application below, $\gamma = \lambda$), it is in the
horizontal direction; it lasts $\gamma \le \lambda$ steps but will be
used in $\le_{K^2_f}$-increasing subsequences of $\langle \bold k_i:i
< \lambda^+\rangle$.
\begin{definition}
\label{it.14}  
For $\delta < \lambda$ and $\gamma \le \lambda$ 
let $\Game_{\delta,\gamma}$ be the following game between the player
INC (incomplete) and COM (complete).

A play last $\gamma$ moves.  In the $\beta$-th move a pair
$(\bold p_\beta,\bold q_\beta)$ is chosen such that $\bold p_\beta
\le^+_{K_1} \bold q_\beta$ and $\beta(1) < \beta \Rightarrow
(\bold p_{\beta(1)} \le_{K_1} \bold p_\beta) \wedge (\bold q_{\beta(1)}
\le_{K_1} \bold q_\beta) \wedge (u_{\bold p_\beta} \cap u_{\bold
q_{\beta(1)}} = u_{\bold p_{\beta(1)}})$ and $u_{\bold p_\beta} \cap
\lambda = \delta$ and $u_{\bold q_\beta} \cap \lambda = u_{\bold q_0}
\cap \lambda \supseteq \delta +1$.

In the $\beta$-th move \underline{first} INC chooses $(\bold
p_\beta,u_\beta)$ such that $\bold p_\beta$ satisfies the requirements
and $u_\beta$ satisfies the requirements on $u_{\bold q_\beta}$
(i.e. $\cup\{u_{\bold q_\alpha}:\alpha < \beta\} \cup
u_{\bold p_\beta} \subseteq u_\beta \in [\lambda^+]^{< \lambda}$ 
and $u_\beta \cap \lambda = u_{\bold q_0} \cap 
\lambda)$ and say $u_\beta \backslash 
u_{\bold p_\beta} \backslash \cup\{u_{\bold q_\gamma}:
\gamma < \beta\}$ has cardinality $\ge |\delta|$
(if $\lambda$ is weakly inaccessible we may be interested in asking more).

\noindent
\underline{Second}, COM chooses $\bold q_\beta$ as required such that 
$u_\beta \subseteq u[\bold q_\beta]$.

A player which has no legal moves loses the play, and arriving to the
$\gamma$-th move, COM wins.
\end{definition}

\begin{remark}
\label{it.16}  
It is not problematic for COM to have a winning
strategy.  But having ``interesting" winning strategies is the crux of
the matter.  More specifically, any
application of this section is by choosing such strategies.

Such examples are the
\mn
\begin{enumerate}
\item[$(a)$]   \underline{lazy man strategy}: preserve 
$\bbP_{\bold q_\beta} = \bbP_{\bold q_0} \times_{\bbP_{\bold p_0}} 
\bbP_{\bold p_\beta}$ recalling Claim \ref{it.8}
\sn
\item[$(b)$]   it is never too late to become lazy, i.e. arriving to
$(\bold p_{\beta(*)},\bold q_{\beta(*)})$ the COM player may decide
that $\beta \ge \beta(*) \Rightarrow \bbP_{\bold q_\beta} =
\bbP_{q_{\beta(*)}} \times_{\bbP_{\bold p_{\beta(*)}}} \bbP_{\bold p_\beta}$
\sn
\item[$(c)$]   definable forcing strategy, i.e. preserve ``$\bbP_{\bold
q_\beta}/\bbP_{\bold p_\beta}$ is a definable c.c.c. forcing (in
$\bold V^{\bbP[\bold p_\beta]}$)".
\end{enumerate}
\end{remark}

\begin{definition}
\label{it.20} 
We say $f$ is $\lambda$-appropriate if
\mn
\begin{enumerate}
\item[$(a)$]   $f \in {}^\lambda(\lambda+1)$
\sn
\item[$(b)$]   $\alpha < \lambda \wedge f(\alpha) < \lambda \Rightarrow 
(\exists \beta)[f(\alpha) = \beta +1]$
\sn
\item[$(c)$]   if $\varepsilon < \lambda^+,\langle u_\alpha:\alpha <
\lambda\rangle$ is an increasing continuous sequence of subsets of
$\varepsilon$ of cardinality $< \lambda$ with union $\varepsilon$
\then \, $\{\delta < \lambda$: otp$(u_\delta) < f(\delta)\}$ is a
stationary subset of $\lambda$.
\end{enumerate}
\end{definition}

\begin{convention}
\label{it.20b}
Below $f$ is $\lambda$-appropriate function.
\end{convention}

\noindent
We arrive to defining the set of approximations of size $\lambda$ (in
the main application $f_*$ is constantly $\lambda$); we shall later
connect it to the oracle version (also see the introduction).
\begin{definition}
\label{it.21}  
For $f_*$ a $\lambda$-appropriate function let 
$K^2_{f_*}$ be the family of $\bold k$ such that:
\mn
\begin{enumerate}
\item[$(a)$]   $\bold k = \langle E,\bar{\bold p},S,\bar{\bold s},
\bar{\bold g},f\rangle$ 
\sn
\item[$(b)$]  $E$ is a club of $\lambda$
\sn
\item[$(c)$]   $\bar{\bold p} = \langle \bold p_\alpha:\alpha \in E\rangle$
\sn
\item[$(d)$]   $\bold p_\alpha \in K^1_\lambda$
\sn
\item[$(e)$]    $\bold p_\alpha \le_{K_1} \bold p_\beta$ for $\alpha
< \beta$ from $E$
\sn
\item[$(f)$]   if $\delta \in \text{ acc}(E)$ then $\bold p_\delta =
\cup\{\bold p_\alpha:\alpha \in E \cap \delta\}$
\sn
\item[$(g)$]   $S \subseteq \lambda$ is a stationary set of limit ordinals
\sn
\item[$(h)$]   if $\delta \in S \cap E$ (hence a limit ordinal) 
then $\delta +1 \in E$
\sn
\item[$(i)$]   $\bar{\bold s} = \langle \bold s_\delta:\delta \in E
\cap S\rangle$
\sn
\item[$(j)$]   $\bold s_\delta$ is a winning strategy for the 
player COM in $\Game_{\delta,f(\delta)}$, see \ref{it.24}(1)
\sn
\item[$(k)$]   $\bar{\bold g} = \langle \bold g_\delta:\delta \in S
\cap E\rangle$
\sn
\item[$(l)$]   $\bullet \quad \bold g_\delta$ is an 
initial segment of a play of $\Game_{\delta,f_*(\delta)}$ in which 
the COM player 

\hskip25pt uses the strategy $\bold s_\delta$
\sn
\item[${{}}$]  $\bullet \quad$ if its length is $< f_*(\delta)$ then
$\bold g_\delta$ has a last move
\sn
\item[${{}}$]  $\bullet \quad (\bold p_\delta,\bold p_{\delta
+1})$ is the pair chosen in the last move, call it mv$(\bold
g_\delta)$
\sn
\item[${{}}$]  $\bullet \quad$ let $S^0 = 
\{\delta \in S \cap E:\bold g_\delta$ has
length $< f_*(\delta)\}$ and $S^1 = S \cap E \backslash S_0$
\sn
\item[$(m)$]   if $\alpha < \beta$ are from $E$ \then  \, $\bold p_\alpha
\le^+_{K_1} \bold p_\beta$, so in particular $\bbP_\beta/\bbP_\alpha$ is
absolutely c.c.c. that is if $\bbP \lessdot \bbP'$ and $\bbP'$
is c.c.c. then $\bbP' *_{{\bbP}_\alpha} \bbP_\beta$ is c.c.c.;
this strengthens clause (e)
\sn
\item[$(n)$]   $f \in {}^\lambda \lambda$
\sn
\item[$(o)$]  if $\delta \in S \cap E$ then $f(\delta)+1$
is the length of $\bold g_\delta$
\sn
\item[$(p)$]   for every $\delta \in E$, if $f_*(\delta) < \lambda$
\then \, $f(\delta) \le \text{ otp}(u_{\bold p_\delta})$.
\end{enumerate}
\end{definition}

\begin{remark}
\label{it.24} 
1) Concerning clause (j), recall (using
the notation of Definition \ref{it.14}) that during a play 
the player INC chooses
$\bold p_\varepsilon$ and COM chooses $\bold q_\varepsilon,\varepsilon
\le f(\delta)$ and recalling clause $(o)$ we see that 
$(\bold p_{f(\delta)},\bold q_{f(\delta)})$ there
stands for $(\bold p_\delta,\bold p_{\delta +1})$ here.  You may
wonder from where does the $(\bold p_\varepsilon,\bold q_\varepsilon)$
for $\varepsilon < f(\delta)$ comes from; the answer is that you
should think of $\bold k$ as a stage in an increasing sequence of
approximations of length $f(\delta)$ and $(\bold p_\varepsilon,\bold
q_\varepsilon)$ comes from the $\delta$-place in the
$\varepsilon$-approximation.  This is cheating a bit - the sequence of
approximations has length $< \lambda^+$, but as on a club of $\lambda$
this reflects to length $< \lambda$, all is O.K.

\noindent
2) Below we define the partial order 
$\le_{K_2}$ (or $\le_{K^2_{f_*}}$) on the set $K^2_{f_*}$, recall 
our goal is to choose an
$\le_{K_2}$-increasing sequence $\langle \bold
k_\varepsilon:\varepsilon < \lambda^+\rangle$ and our final forcing
will be $\cup\{\Bbb P_{\bold k_\varepsilon}:\varepsilon < \lambda^+\}$.

\noindent
3) Why clause (d) in Definition \ref{it.28}(2) below?  It
is used in the proof of the limit existence claim \ref{it.49}.  This
is because the club $E_{\bold k}$ may decrease (when increasing 
$\bold k$).

Note that we use $\le^*_{K^1_f}$ ``economically".  We cannot in
general demand (in \ref{it.28}(2) below) that for $\alpha < \beta$ from
$E_{\bold k_2} \backslash \alpha(*)$ we have $(\bold p^{\bold
k_1}_\alpha,\bold p^{\bold k_1}_\beta) \le^*_{K_1} (\bold p^{\bold
k_2}_\alpha,\bold p^{\bold k_2}_\beta)$ as the strategies $\bold
s_\delta$ may defeat this.  How will it still help?  Assume $\langle
\bold k_\varepsilon:\varepsilon < \varepsilon(*)\rangle$ is
increasing, $\varepsilon(*) < \lambda$ for simplicity and $\gamma \in
\cap\{E_{\bold k_\varepsilon}:\varepsilon < \varepsilon(*)\} \cap \bigcap
\{S_{\bold k_\varepsilon}:\varepsilon < \varepsilon(*)\} \backslash
\cup\{\alpha(\bold k_\varepsilon,\bold k_\zeta):\varepsilon < \zeta <
\varepsilon(*)\}$ and $\gamma_\varepsilon = \text{ Min}(E_{\bold
k_\varepsilon} \backslash (\gamma +1))$ for $\varepsilon <
\varepsilon(*)$.  We shall have $\langle
\gamma_\varepsilon:\varepsilon < \varepsilon(*)\rangle$ is increasing;
there may be $\delta \in
(\gamma_\varepsilon,\gamma_{\varepsilon +1})$ where $\bold s_\delta$
was active between $\bold k_\varepsilon$ and $\bold k_{\varepsilon
+1}$ but it contributes to $\bbP^{\bold k_{\varepsilon
+1}}_{\gamma_{\varepsilon +1}}/
\bbP^{\bold k_\varepsilon}_{\gamma_\varepsilon}$. 

\noindent
4) If we omit the restriction $u \in [\lambda^+]^{< \lambda}$ and
$f:\lambda \rightarrow \delta^* +1$, replace the club $E$ by an end
segment, we can deal with sequences of length $\delta^*$.

In the direct order in \ref{it.28}(3) we have
$\alpha(*)=0$.  Using e.g. a stationary non-reflecting $S \subseteq
S^{\delta^*}_\lambda$ we can often allow $\alpha(*) \ne 0$.

\noindent
5) Is the ``$\bold s_\delta$ a winning strategy" in addition for
telling us what to do, crucial?  The point is preservation of
c.c.c. in limit of cofinality $\aleph_1$.

\noindent
6) If we use $f_* \in {}^\lambda(\lambda +1)$ constantly $\lambda$, we
do not need $f_{\bold k}$ so we can omit clauses (n),(o),(p) of
\ref{it.21} and (c), and part of another in \ref{it.28}.

\noindent
6A) Alternatively we can omit clause $(o)$ in \ref{it.21} but demand
``$\prod\limits_{\alpha < \lambda} f(\alpha)/\cD$ is
$\lambda^+$-directed", fixing a normal filter $\cD$ on $\lambda$ (and
demand $S_{\bold k} \in \cD^+$).

\noindent
7) The ``omitting type" argument here comes from using the strategies.
\end{remark}

\begin{definition}
\label{it.28}  
1) In Definition \ref{it.21}, let $E = E_{\bold k},
\bar{\bold p} = \bar{\bold p}_{\bold k},
\bold p_\alpha = \bold p^{\bold k}_\alpha = \bold p_\alpha[\bold k],
\bbP_\alpha = \bbP^{\bold k}_\alpha = \bbP_{\bold p_\alpha[\bold k]},S
= S_{\bold k},S_{[\ell]} = S_{\bold k,\ell}$ for $\ell=0,1$, etc. and
we let $\bbP_{\bold k} = \cup\{\bbP^{\bold k}_\alpha:
\alpha \in E_{\bold k}\}$ and $u_{\bold k} = u[\bold k] = 
\cup\{u_{\bold p^{\bold k}_\alpha}:\alpha \in E_{\bold k}\}$.

\noindent
2) We define a two-place relation $\le_{K^2_f}$ on $K^2_f:\bold k_1
 \le_{K^2_f} \bold k_2$ \underline{iff} (both are from $K^2_f$ and)
for some $\alpha(*) < \lambda$
(and $\alpha(\bold k_1,\bold k_2)$ is the first such $\alpha(*)$) we have:
\mn
\begin{enumerate}
\item[$(a)$]   $E_{\bold k_2} \backslash E_{\bold k_1}$ is bounded
 in $\lambda$, moreover $\subseteq \alpha(*)$ 
\sn
\item[$(b)$]    for $\alpha \in E_{\bold k_2} \backslash \alpha(*)$
we have $\bold p^{\bold k_1}_\alpha \le_{K_1} \bold p^{\bold k_2}_\alpha$
\sn
\item[$(c)$]   if $\alpha \in E_{\bold k_2} \backslash \alpha(*)$ then
$f_{\bold k_1}(\alpha) \le f_{\bold k_2}(\alpha)$
\sn
\item[$(d)$]   if $\gamma_0 < \gamma_1 \le \gamma_2 < \lambda,\gamma_0
\in E_{\bold k_2} \backslash (\alpha(*) \cup S_{\bold k_1}),
\gamma_1 = \text{ min}(E_{\bold k_1} \backslash (\gamma_0 +1))$ and
$\gamma_2 = \text{ min}(E_{\bold k_2} \backslash (\gamma_0 +1))$,
\then \, $(\bold p^{\bold k_1}_{\gamma_0},
\bold p^{\bold k_2}_{\gamma_0}) \le_{K_1} (\bold p^{\bold k_1}_{\gamma_1},\bold
p^{\bold k_2}_{\gamma_2})$, see Definition \ref{it.4}(1)
\sn
\item[$(e)$]  if $\delta \in S_{\bold k_1} \cap E_{\bold k_2} \backslash
\alpha(*)$ \then \, $\delta \in S_{\bold k_2} \cap E_{\bold
k_2} \backslash \alpha(*)$; but note that if $f_{\bold k_1}(\delta)
\ge f(\delta)$ we put $\delta$ into $S_{\bold k_2}$ just for
notational convenience
\sn
\item[$(f)$]   if $\delta \in S_{\bold k_1} \cap E_{\bold k_2}
\backslash \alpha(*)$ \then \, $\bold s^{\bold k_2}_\delta = \bold
s^{\bold k_1}_\delta$ and $\bold g^{\bold k_1}_\delta$ is an initial
segment of $\bold g^{\bold k_2}_\delta$
\sn
\item[$(g)$]   if ${\bold k_1} \ne {\bold k_2}$ \then \, 
$u[{\bold k_1}] \ne u[\bold k_2]$
\sn
\item[$(h)$]   if $\alpha < \beta$ are from 
$E_{\bold k_2} \backslash \alpha(*)$ \then \,
$(\bold p^{\bold k_1}_\alpha,\bold p^{\bold k_2}_\alpha) \le'_{K_1} 
(\bold p^{\bold k_1}_\beta,\bold p^{\bold k_2}_\beta)$, see
Definition \ref{it.4}(2), i.e. if $p \in 
\bbP_{\bold p_\alpha[\bold k_1]},q \in \bbP_{\bold p_\alpha[\bold k_2]}$ and $p
\Vdash_{\bbP_{\bold p_\alpha[\bold k_1]}} ``q \in \bbP_{\bold
p_\alpha[\bold k_2]} /\name G_{\bbP_{\bold p_\alpha[\bold k_1]}}"$
\then \, $p \Vdash_{\bbP_{\bold p_\beta[\bold k_1]}} ``q \in 
\bbP_{\bold p_\beta[\bold k_2]} /\name G_{\bbP_{\bold p_\beta[\bold k_1]}}"$.
\end{enumerate}
\mn
3) We define a two-place relation $\le^{\text{dir}}_{K^2_f}$ on
$K^2_f$ as follows: $\bold k_1 \le^{\text{dir}}_{K^2_f} \bold k_2$ iff
\mn
\begin{enumerate}
\item[$(a)$]   $\bold k_1 \le_{K^2_f} \bold k_2$
\sn
\item[$(b)$]   $E_{\bold k_2} \subseteq E_{\bold k_1}$; no real harm
here if we add $\bold k_1 \ne \bold k_2 \Rightarrow 
E_{\bold k_2} \subseteq \text{ acc}(E_{\bold k_1})$
\sn
\item[$(c)$]   $\alpha(\bold k_1,\bold k_2) = \text{ Min}(E_{\bold k_2})$.
\end{enumerate}
\mn
4) We write
 $K^2_\lambda,\le_{K^2_\lambda},\le^{\text{dir}}_{K^2_\lambda}$ or
 just $K_2,\le_{K_2},<^{\text{dir}}_{K_2}$ 
for $K^2_f,\le_{K^2_f},\le^{\text{dir}}_{K^2_f}$ 
when $f$ is constantly $\lambda$.
\end{definition}

\begin{remark}
\label{it.30}
1) In \cite{Sh:669} we may increase $S$ as well as here but we may
replace clause (e) by
\mn
\begin{enumerate}
\item[$(e)'$]   $\delta \in S_{\bold k_1} \cap E_{\bold k_2} \backslash
\alpha(*)$ \underline{iff} $f_{\bold k_2}(\delta) < f(\delta) \wedge \delta
\in S_{\bold k_2} \cap E_{\bold k_2} \backslash \alpha(*)$.
\end{enumerate}
\mn
If we do this, is it a great loss?  No!
This can still be done here by choosing $\bold s_\delta$ 
such that as long as INC chooses
$u_\beta$ of certain form (e.g. $u_\beta \backslash u^{\bold p_\beta}
= \{\delta\}$) the player COM chooses $\bold q_\beta = \bold p_\beta$.
If we can allow in Definition \ref{it.28}(2) to extend $S$ but a priori 
start with $\langle S_\varepsilon:\varepsilon < \lambda^+\rangle$ 
such that $S_\varepsilon
\subseteq \lambda$ and $S_\varepsilon \backslash S_\zeta$ is bounded
in $\lambda$ when $\varepsilon < \zeta < \lambda$.

\noindent
2) We can weaken clause $(e)$ of \ref{it.28}(2) to
\mn
\begin{enumerate}
\item[$(e)''$]  if $\delta \in S_{\bold k_1} \cap E_{\bold k_2}
   \backslash \alpha(*)$ and $f_{\bold k_2}(\delta) < f(\delta)$ then
   $\delta \in S_{\bold k_2}$.
\end{enumerate}
\mn
But then we have to change accordingly, e.g. \ref{it.28}(c),(f),
\ref{it.35}(c).

\noindent
3) We can define $\bold k_1 \le_{K^2_f} \bold k_2$ demanding
$(S_{\bold k_1},\bar{\bold s}_{\bold k_1}) = (S_{\bold k_2},
\bar{\bold s}_{\bold k_2})$ but replace everywhere ``$\delta \in S_{\bold k}
\cap E_{\bold k}"$ by ``$\delta \in S_{\bold k} \cap E_{\bold k}
\wedge f_{\bold k}(\delta) \le f(\delta)$" so omit clause (e) of \ref{it.28}.
\end{remark}

\begin{observation}
\label{it.31}  
1) $\le_{K^2_f}$ is a partial order on $K^2_f$.

\noindent
2) $\le^{\text{\rm dir}}_{K^2_f} \subseteq \le_{K^2_f}$ is a partial order
on $K^2_f$.

\noindent
3) If $\bold k_1 \le_{K^2_f} \bold k_2$ then $\bbP_{\bold k_1} \lessdot
\bbP_{\bold k_2}$.

\noindent
4) If $(\bold k_\varepsilon:\varepsilon < \lambda^+\rangle$ is
$<_{K^2_f}$-increasing and $\bbP = \cup\{\bbP_{\bold
k_\varepsilon}:\varepsilon < \lambda^+\}$ \then
\mn
\begin{enumerate}
\item[$(a)$]   $\bbP$ is a c.c.c. forcing notion of cardinality
$\le \lambda^+$
\sn
\item[$(b)$]   $\bbP_{\bold k_\varepsilon} \lessdot \bbP$ for
$\varepsilon < \lambda^+$.
\end{enumerate}
\end{observation}

\begin{definition}
\label{it.35}  
1) Assume $\bar{\bold k} = \langle \bold k_\varepsilon:
\varepsilon < \varepsilon(*)\rangle$ is
$\le_{K^2_f}$-increasing with $\varepsilon(*)$ a limit ordinal $<
\lambda$.  We say $\bold k$ is a limit of $\bar{\bold k}$ \when \,
$\varepsilon < \varepsilon(*) \Rightarrow \bold k_\varepsilon
\le_{K^2_f} \bold k \in K^2_f$ and for some $\alpha(*)$
\mn
\begin{enumerate}
\item[$(a)$]   $\alpha(*) = \cup\{\alpha(\bold k_\varepsilon,\bold
 k_\zeta):\varepsilon < \zeta < \varepsilon(*)\}$
\sn
\item[$(b)$]    $E_{\bold k} \backslash \alpha(*) \subseteq \cap\{E_{\bold
k_\varepsilon} \backslash \alpha(*):\varepsilon < \varepsilon(*)\}$
\sn
\item[$(c)$]  $S_{\bold k} = (\cup\{S_{\bold k_\varepsilon}:
\varepsilon < \varepsilon(*)\}) \cap (\cap\{E_{\bold k_\varepsilon}:
\varepsilon < \varepsilon(*)\}) \backslash \alpha(*)$
\sn
\item[$(d)$]   if $\delta \in S_{\bold k}$ \then \,
$\bold g^{\bold k_\varepsilon}_\delta$ is an initial segment of
$\bold g^{\bold k}_\delta$ for every $\varepsilon < \varepsilon(*)$
\sn
\item[$(e)$]   $f_{\bold k}(\delta) = \cup\{f_{\bold
k_\varepsilon}(\delta):\varepsilon < \varepsilon(*)\} + 1$ for $\delta
\in S_{\bold k}$.
\end{enumerate}
\mn
2) Assume $\bar{\bold k} = \langle \bold k_\varepsilon:
\varepsilon < \lambda\rangle$ is
$\le_{K^2_f}$-increasing continuous.  
We say $\bold k$ is a limit of $\bar{\bold k}$ \when \, 
$\varepsilon < \lambda \Rightarrow \bold k_\varepsilon
\le \bold k \in K^2_f$ and for some $\bar \alpha$
\mn
\begin{enumerate}
\item[$(a)$]   $\bar \alpha = \langle \alpha_\varepsilon:\varepsilon
< \lambda\rangle$ is increasing continuous, $\lambda >
\alpha_\varepsilon \in \cap\{E_{\bold k_\zeta}:\zeta < \varepsilon\}
\backslash \cup\{\alpha(\bold k_{\zeta_1},\bold k_{\zeta_2}):\zeta_1 < \zeta_2
< 1 +  \varepsilon\}$
\sn
\item[$(b)$]    $E_{\bold k} = \{\alpha_\varepsilon:\varepsilon <
\lambda\} \cup \{\alpha_\varepsilon +1:\varepsilon < \lambda$ and
$\varepsilon \in S\}$
\sn
\item[$(c)$]   $S_{\bold k} = \{\alpha_\varepsilon:\alpha_\varepsilon 
\in \bold S_{\bold k_\zeta}$ for every $\zeta < \varepsilon$ large enough$\}$
\sn
\item[$(d)$]  if $\delta = \alpha_\varepsilon \in S_{\bold
k_\varepsilon}$ then $\bold g^{\bold k}_\delta = 
\bold g^{\bold k_\varepsilon}_\delta$
\sn
\item[$(e)$]   if $\alpha < \delta$ and $\zeta = \text{
Min}\{\varepsilon:\alpha \le \alpha_{\varepsilon +1}\}$ 
then $f_{\bold k}(\alpha) = f_{\bold k_\zeta}(\alpha)$.
\end{enumerate}
\mn
3) We say that $\langle \bold k_\varepsilon:\varepsilon <
\varepsilon(*)\rangle$ is $\le_{K^2_f}$-increasing continuous
\when \, :
\mn
\begin{enumerate}
\item[$(a)$]   $\bold k_\varepsilon \le_{K^2_f} \bold k_\zeta$ for
 $\varepsilon < \zeta < \varepsilon(*)$
\sn
\item[$(b)$]   $\bold k_\varepsilon$ is a limit of $\langle \bold
k_{\xi(\zeta)}:\zeta < \text{ cf}(\varepsilon)\rangle$ for some
increasing continuous sequence $\langle \xi(\zeta):\zeta < \text{
cf}(\varepsilon)\rangle$ of ordinals with limit $\varepsilon$, for every
limit $\varepsilon < \varepsilon(*)$, by part (1) or part (3).
\end{enumerate}
\end{definition}

\begin{definition}
\label{it.38}  
1) In part (1) of \ref{it.35},
we say ``a direct limit" when in addition
\mn
\begin{enumerate}
\item[$(\alpha)$]   the sequences are $\le^{\text{dir}}_{K^2_f}$-increasing
\sn
\item[$(\beta)$]    in clause $(a),(b)$ we have equality
\sn
\item[$(\gamma)$]   $\bold p^{\bold k}_{\text{min}(E_{\bold k})}$ is
the exact union of $\langle \bold p^{\bold k}_{\text{min}(E_{\bold
k_\varepsilon})}:\varepsilon < \varepsilon(*)\rangle$
\sn
\item[$(\delta)$]   if $\gamma \in E_{\bold k},\xi <
\varepsilon(*),\gamma \notin S^1_{\bold k_\xi}$
and $\langle \gamma_\varepsilon:\varepsilon \in
[\xi,\varepsilon(*)]\rangle$ is defined by $\gamma_\xi =
\xi,\gamma_\varepsilon = \text{ min}(E_{\bold k_\varepsilon}
\backslash (\gamma +1))$, so $\langle \gamma_\varepsilon:\varepsilon
\in [\xi,\varepsilon(*)]\rangle$ is an $\le$-increasing continuous sequence
of ordinals, \then \, $\bold p^{\bold k}_{\gamma_{\varepsilon(*)}}/
\bold p^{\bold k}_\gamma = \cup\{\bold p^{\bold
k}_{\gamma_\varepsilon}/\bold p^{\bold k}_\gamma:\varepsilon
\in [\xi,\varepsilon(*))\}$ with the obvious meaning. 
\end{enumerate}
\mn
2) In part (2) of Definition \ref{it.35} we say 
a ``direct limit" when in addition
\mn
\begin{enumerate}
\item[$(\alpha)$]   the sequence is $\le^{\text{dir}}_{K^2_f}$
\sn
\item[$(\beta)$]  $\alpha_\varepsilon = 
\text{ Min}(E_{\bold k_\varepsilon})$ \underline{or} the
$\varepsilon$-th member of $E_{\bold k_\varepsilon}$.
\end{enumerate}
\mn
3) We say that $\bar{\bold k} = \langle \bold k_\varepsilon:\varepsilon <
\varepsilon(*)\rangle$ is $\le^{\text{dir}}_{K^2_f}$-increasing
continuous or directly increasing continuous \when \, :
\mn
\begin{enumerate}
\item[$(a)$]   $\bold k_\varepsilon \le^{\text{dir}}_{K^2_f} 
\bold k_\zeta$ for $\varepsilon \le \zeta < \varepsilon(*)$
\sn
\item[$(b)$]   if $\varepsilon < \varepsilon(*)$ is a limit ordinal
then $\bold k_\varepsilon$ is a (really the) direct limit of
$\bar{\bold k} \restriction \varepsilon$.
\end{enumerate}
\end{definition}

\begin{claim}
\label{it.42}  If $\bold k_1 \le_{K^2_f} \bold k_2$
\then \, for some $\bold k'_2$ we have
\mn
\begin{enumerate}
\item[$(a)$]   $\bold k_1 \le^{\text{\rm dir}}_{K^2_f} \bold k'_2$ 
\sn
\item[$(b)$]  $\bold k_2 \le_{K^2_f} \bold k'_2 \le_{K^2_f}
\bold k_2$
\sn
\item[$(c)$]   $\bold k_2,\bold k'_2$ are almost equal - the only
differences being $E_{\bold k'_2} = E_{\bold k_2} \backslash 
\text{\rm min}(E_{\bold k'_2}),S_{\bold k'_2} \subseteq S_{\bold
k_2}$, etc.
\end{enumerate}
\end{claim}
\bigskip

\begin{claim}
\label{it.49}
\underline{The limit existence claim}  
1) If $\varepsilon(*) < \lambda$ is a limit ordinal and $\bar{\bold k} =
\langle \bold k_\varepsilon:\varepsilon < \varepsilon(*)\rangle$ is a
[directly] increasing continuous \then \, $\bar{\bold k}$ has a 
[direct] limit.

\noindent
2) Similarly for $\varepsilon(*) = \lambda$.
\end{claim}

\begin{PROOF}{\ref{it.49}}
It is enough to prove the direct version.

\noindent
1) We define $\bold k = \bold k_{\varepsilon(*)}$ as in
the definition, we have no freedom left.

The main points concern the c.c.c. and the absolute c.c.c.,
$\le'_{K^0_1},\le_{K_1}$ demands.  
We prove the relevant demands by induction on $\beta \in 
E_{\bold k_{\varepsilon(*)}}$.
\medskip

\noindent
\underline{Case 1}:  $\beta = \text{ min}(E_{\bold k_{\varepsilon(*)}})$.

First note that $\langle \bold p^\varepsilon_{\text{min}(E_{\bold
k_\varepsilon})}:\varepsilon \le \varepsilon(*)\rangle$ is increasing
continuous (in $K^1_\lambda$) moreover $\langle 
\bbP[\bold p^{\bold k_\varepsilon}_{\text{min}(E_{\bold
k_\varepsilon})}]:\varepsilon \le
\varepsilon(*)\rangle$ is increasing continuous, see clause $(\gamma)$
of Definition \ref{it.38}(1).  As each 
$\bbP[\bold p_{\text{min}(E_{\bold k_\varepsilon)}}]$ is c.c.c. if
$\varepsilon < \varepsilon(*)$, we know that this holds for
$\varepsilon = \varepsilon(*)$, too.
\medskip

\noindent
\underline{Case 2}:  $\beta = \delta +1,
\delta \in S^1_{\bold k} \cap E_{\bold k}$.

Since $\bold s^{\bold k}_\delta$ is a winning strategy in
the game $\Game_{\delta,f(\delta)}$ we have 
$\bold p^{\bold k_{\varepsilon(*)}}_\delta \le^+_{K_1} 
\bold p^{\bold k_{\varepsilon(*)}}_\beta$.  
But what if the play is over?  Recall that
in Definition \ref{it.20},
$f(\delta) = \lambda$ or $f(\delta)$ is successor and $\langle
f_{\bold k_\varepsilon}(\delta):\varepsilon < \varepsilon(*)\rangle$ is
(strictly) increasing, so this never happens; it may happen when we
try to choose $\bold k'$ such that $\bold k <_{K^2_f} \bold k'$, see
\ref{it.56}. 

We also have to show: if $\alpha \in \beta \cap E_{\bold k}$ 
then $\bbP[\bold p^{\bold k}_\beta]/\bbP[\bold p^{\bold k}_\alpha]$ is
absolutely c.c.c.  First, if $\alpha = \delta$
this holds by Definition \ref{it.7}(3) of $\le^+_{K_1}$ and the
demand $\bold p_\beta \le^+_{K_1} \bold q_\beta$ in Definition
\ref{it.14} (and clause $(\ell)$ of Definition \ref{it.21}).  
Second, if $\alpha < \delta$, it is enough to show
that $\bbP[\bold p^{\bold k}_\beta]/\bbP[\bold p^{\bold k}_\delta]$ and
$\bbP[\bold p^{\bold k}_\delta]/\bbP[\bold p^{\bold k}_\alpha]$ are
absolutely c.c.c., but the first holds by the previous sentence, the
second by the induction hypothesis.
In particular, when $\varepsilon < \varepsilon(*) \Rightarrow 
\bbP^{\bold k_\varepsilon}_\beta \lessdot \bbP^{\bold k}_\beta$.
\medskip

\noindent
\underline{Case 3}:  For some 
$\gamma,\gamma = \text{ max}(E_{\bold k} \cap
\beta),\gamma \notin S^1_{\bold k}$.

As $\gamma \notin S_{\bold k}$ there is $\xi < \varepsilon(*)$
such that $\gamma \notin S^1_{\bold k_\xi}$ let $\gamma_\xi = \gamma$
and for $\varepsilon \in (\xi,\varepsilon(*)]$ we define
$\gamma_\varepsilon =: \text{ min}(E_{\bold k_\varepsilon} \backslash
(\beta +1))$.  Now as $\bar{\bold k}$ is directly increasing 
continuous we have
\mn
\begin{enumerate}
\item[$\circledast$]  $(a) \quad \langle
\gamma_\varepsilon:\varepsilon \in [\xi,\varepsilon(*)]\rangle$ is
increasing continuous
\sn
\item[${{}}$]  $(b) \quad \gamma_\xi = \gamma$
\sn
\item[${{}}$]   $(c) \quad \gamma_{\varepsilon(*)} = \beta$
\sn
\item[${{}}$]   $(d) \quad \langle \bold p^{\bold
k_\varepsilon}_{\gamma_\varepsilon}:\varepsilon \in
[\xi,\varepsilon(*)]\rangle$ is increasing continuous.
\end{enumerate}
\mn
So by claim \ref{it.5} we are done, the main point is that clause (d)
there holds by clause (d) of the definition of $\le_{K^2_f}$ in
\ref{it.28}(2).
\medskip

\noindent
\underline{Case 4}:  $\beta = \sup(E_{\bold k} \cap \beta)$.

It follows by the induction hypothesis and \ref{it.3}(3) as $\langle
\bold p^{\bold k}_\gamma:\gamma \in E_k \cap \beta\rangle$ is
$\le^+_{K_1}$-increasing continuous with union 
$\bold p^{\bold k}_\beta$; of course we use clause (h) of Definition
\ref{it.28}, so Definition \ref{it.4}(2),(5) applies.  

\noindent
2) Similarly. 
\end{PROOF}

\noindent
The following is an atomic step toward having MA$_{< \lambda}$.
\begin{claim}
\label{it.56}  
Assume
\mn
\begin{enumerate}
\item[$(a)$]   $\bold k_1 \in K^2_f$
\sn
\item[$(b)$]  $\alpha(*) \in E_{\bold k_1}$
\sn
\item[$(c)$]  $\name{\bbQ}$ is a 
$\bbP[\bold p^{\bold k_1}_{\alpha(*)}]$-name of a c.c.c. 
forcing (hence $\Vdash_{\bbP_{\bold k_1}} ``\name{\bbQ}$ is a 
c.c.c. forcing")
\sn
\item[$(d)$]   $u_* \subseteq \lambda^+$ is disjoint to $u[\bold k_1]
= \cup\{u_{\bold p_\alpha[\bold k_1]}:\alpha \in E_{\bold k}\}$ and of
cardinality $< \lambda$ but $\ge |\name{\bbQ}|$.
\end{enumerate}
\mn
\Then \, we can find $\bold k_2$ such that
\mn
\begin{enumerate}
\item[$(\alpha)$]   $\bold k_1 \le^{\text{dir}}_{K^2_f} 
\bold k_2 \in K^2_f$
\sn
\item[$(\beta)$]  $E_{\bold k_2} = E_{\bold k_1} \backslash \alpha(*)$
\sn
\item[$(\gamma)$]  $u^{\bold k_2}_\alpha = u^{\bold k_1}_\alpha
\cup u_*$ for $\alpha \in E_{\bold k_2} \cap S^1_{\bold k_1}$
\sn
\item[$(\delta)$]  $\bbP_{\bold p_{\alpha(*)}[\bold k_2]}$ is
isomorphic to $\bbP_{\bold p_{\alpha(*)}[\bold k_1]}  *
\name{\bbQ}$ over $\bbP_{\bold p_{\alpha(*)}[\bold k_1]}$
\sn
\item[$(\varepsilon)$]  $S_{\bold k_2} = S_{\bold k_1}$ and 
$\bar{\bold s}_{\bold k_2} = \bar{\bold s}_{\bold k_1} \rest 
G_{\bold k_2}$
\sn
\item[$(\zeta)$]  $f_{\bold k_2} = f_{\bold k_1} +1$
\sn
\item[$(\eta)$]   if $\Vdash_{\bbP_{k_1} * \name{\bbQ}} 
``\name \rho \in {}^\omega 2$ but $\name \rho \notin 
\bold V[\name G_{\bbP_{\bold k_1}}]"$ then $\Vdash_{\bbP_{\bold k_2}}
``\name \rho \in {}^\omega 2$ but $\name \rho \notin 
\bold V[\name G_{\bbP_{\bold k_1}}]$ provided that the 
strategies preserve this which they do under the criterion here.
\end{enumerate}
\end{claim}

\begin{PROOF}{\ref{it.56}}
We choose $\bold p^{\bold k_2}_\alpha$ by induction on
$\alpha \in E_{\bold k_1} \backslash \alpha(*)$, keeping all relevant
demands (in particular $u_{\bold p_\alpha[\bold k_2]} \cap u[\bold
k_1] = u_{\bold p_\alpha[\bold k_1]})$.
\bigskip

\noindent
\underline{Case 1}:  $\alpha = \alpha(*)$.

As only the isomorphism type of $\name{\bbQ}$ is
important, \wilog \, 
$\Vdash_{\bbP[\bold p^{\bold k_1}_{\alpha(*)}]}$ 
``every member of $\name{\bbQ}$ belongs to $u_*$".

So we can interpret the set of elements of 
$\bbP_{\bold p_{\alpha(*)}[\bold k_1]} * \name{\bbQ}$ such
that it is $\subseteq {\cH}_{< \aleph_1}(u_{\bold p_{\alpha(*)}
[\bold k_1]} \cup u_*)$.

Now $\bbP_{\bold p_{\alpha(*)}[\bold k_1]} \lessdot \bbP_{\bold
p_{\alpha(*)}[\bold k_2]}$ by the classical claims on composition of
forcing notions.
\bigskip

\noindent
\underline{Case 2}:  $\alpha = \delta +1,\delta \in S_{\bold k_1} \cap
E_{\bold k_1} \backslash \alpha(*)$.

The case split to two subcases.
\bigskip

\noindent
\underline{Subcase 2A}:  The play 
$\bold g^{\bold k_1}_\delta$ is not over,
i.e. $f(\delta)$ is larger than the length of the play so far.

In this case do as in case 2 in the proof of \ref{it.49}, just use
$\bold s_\delta$. 
\bigskip

\noindent
\underline{Subcase 2B}:  The play $\bold g^{\bold k_1}_\delta$ is over.

In this case let $\bbP^{\bold k_2}_{\delta +1} = \bbP^{\bold
k_1}_{\delta +1} *_{\bbP^{\bold k_1}_\delta} \bbP^{\bold
k_2}_\delta$, in fact, $\bold p^{\bold k_2}_{\delta +1} = \bold p^{\bold
k_1}_{\delta +1} *_{\bold p^{\bold k_1}_\delta} \bold p^{\bold
k_2}_\delta$ (and choose $u_{\bold p_{\delta +1}[\bold k_2]}$
appropriately).  Now possible and $(\bold p^{\bold k_1}_\delta,\bold
p^{\bold k_2}_\delta) <'_{K_1} (\bold p^{\bold k_1}_{\delta +1},\bold
p^{\bold k_2}_{\delta +1})$ by \ref{it.8}.
\bigskip

\noindent
\underline{Case 3}:  For some $\gamma,\gamma = 
\text{ max}(E_{\bold k} \cap \beta) \ge \alpha(*)$ and 
$\gamma \notin S_{\bold k}$.

Act as in Subcase 2B of the proof of \ref{it.49}
\bigskip

\noindent
\underline{Case 4}:  $\beta = \sup(E_{\bold k} \cap \beta)$.

As in Case 4 in the proof of \ref{it.49}.  
\end{PROOF}
\newpage

\section {${\frak p} = {\frak t}$ does not decide the existence of a
peculiar cut} 

We deal here with a problem raised in \cite{Sh:885}, toward this we
quote from there.  Recall (Definition \cite[1.10]{Sh:885}).
\begin{definition}
\label{pt.0n}  
Let $\kappa_1,\kappa_2$ be infinite regular cardinals.  
A $(\kappa_1,\kappa_2)$-peculiar cut in 
${}^\omega \omega$ is a pair $(\langle f_i:i < \kappa_1\rangle,\langle
f^\alpha:\alpha < \kappa_2\rangle)$ of sequences of functions in
${}^\omega \omega$ such that:
\mn
\begin{enumerate}
\item[$(\alpha)$]  $(\forall i<j<\kappa_1)(f_j 
<_{J^{\text{bd}}_\omega} f_i)$,
\sn
\item[$(\beta)$]  $(\forall \alpha < \beta < \kappa_2)(f^\alpha 
<_{J^{\text{bd}}_\omega} f^\beta)$,
\sn
\item[$(\gamma)$]  $(\forall i < \kappa_1)(\forall \alpha <
\kappa_2)(f^\alpha <_{J^{\text{bd}}_\omega} f_i)$,
\sn
\item[$(\delta)$]  if $f:\omega \rightarrow \omega$ is such that
$(\forall i < \kappa_1)(f \le_{J^{\text{bd}}_\omega} f_i)$, 
then $f \le_{J^{\text{bd}}_\omega} f^\alpha$ for some 
$\alpha < \kappa_2$,
\sn
\item[$(\varepsilon)$]  if $f:\omega \rightarrow \omega$ is such
that $(\forall \alpha < \kappa_2)(f^\alpha \le_{J^{\text{bd}}_\omega}  
f)$, then $f_i \le_{J^{\text{bd}}_\omega} f$ for some $i < \kappa_1$.
\end{enumerate}
\mn
Recall that if $\gp < \gt$ then for some regular $\kappa < \gp$ there
is a $(\kappa,\gp)$-peculiar cut, (\cite[1.12]{Sh:885}).  Also $\gp =
\aleph_1 \Rightarrow \gt = \gp$ by the classicl theorem of Rothenberg
and MA$_{\aleph_1} + \gp = \aleph_2 \Rightarrow \gt = \aleph_2$ by
\cite[2.3]{Sh:885}. 
\end{definition}

\noindent
Recall (from \cite{Sh:885}) that
\begin{claim}
\label{pt.Op}
1) If there is a $(\kappa_1,\kappa_2)$-peculiar \then \, recall from there
   that the motivation of looking at $(\kappa_1,\kappa_2)$-peculiar
   type is understanding the case $\gp > \gt$.

\noindent
1A) In particular, if $\gp < \gt$ \then \, there 
is a $(\kappa_1,\kappa_2)$-peculiar type for some
(regular) $\kappa_1,\kappa_2$ satisfying $\kappa_1 < \kappa_2 =
   \gt$, see \cite{Sh:885}, $\gt \le \gp \le 
\text{\rm max}\{\kappa_1,\kappa_2\}$.

\noindent
2) There is a $(\kappa_1,\kappa_2)$-peculiar cut iff there is a
   $(\kappa_1,\kappa_1)$-peculiar cut.  
\end{claim}

\begin{PROOF}{\ref{pt.OP}}
1) Straight.

\noindent
2) Trivial.
\end{PROOF}

\begin{observation}
\label{pt.0q}   If $(\bar \eta^{\text{up}},
\bar \eta^{\text{dn}})$ is a peculiar $(\kappa_{\text{up}},
\kappa_{\text{dn}})$-cut and if $A \subseteq
\omega$ is infinite, $\eta \in {}^\omega \omega$ \then \,:
\mn
\begin{enumerate}
\item[$(a)$]   $\eta <_{J^{\text{bd}}_A} \eta^{\text{up}}_\alpha$
for every $\alpha < \kappa_{\text{up}}$ \underline{iff} $\eta
<_{J^{\text{bd}}_A} \eta^{\text{dn}}_\beta$ for every large enough
$\beta < \kappa_{\text{dn}}$
\sn
\item[$(b)$]   $\neg(\eta^{\text{up}}_\alpha <_{J^{\text{bd}}_A} \eta)$
for every $\alpha < \kappa_{\text{up}}$ \underline{iff} 
$\neg(\eta^{\text{dn}}_\beta
<_{J^{\text{bd}}_A} \eta)$ for every large enough $\beta < 
\kappa_{\text{dn}}$.
\end{enumerate}
\end{observation}

\begin{PROOF}{\ref{pt.0q}}
\underline{Clause (a)}:   The implication $\Leftarrow$ is trivial
as $\beta <\kappa_{\text{dn}} \wedge \alpha < \kappa_{\text{up}}
\Rightarrow \eta^{\text{dn}}_\beta <_{J^{\text{bd}}_\omega}
\eta^{\text{up}}_\alpha$.  So assume the leftside.

We define $\eta' \in {}^\omega \omega$ by: $\eta'(n)$ is $\eta(n)$ if
$n \in A$ and is $0$ if $n \in \omega \backslash A$.  Clearly $\eta'
<_{J^{\text{bd}}_\omega} \eta^{\text{up}}_\alpha$ for every $\alpha <
\kappa_{\text{up}}$ hence by clause $(\varepsilon)$ of \ref{pt.0n} we
have  $\eta' \le_{J^{\text{bd}}_\omega}
\eta^{\text{dn}}_\beta$ for every large enough $\beta <
\kappa_{\text{dn}}$ hence $\eta = \eta' \restriction A
\le_{J^{\text{bd}}_A} \eta^{\text{dn}}_{\beta +1} 
<_{J^{\text{bd}}_A} \eta^{\text{dn}}_\beta$ for every large enough
$\beta < \kappa_{\text{dn}}$.
\medskip

\noindent
\underline{Clause (b)}:  Again the direction 
$\Leftarrow$ is obvious.  For the other direction
define $\eta' \in {}^\omega \omega$ by $\eta'(n)$ is $\eta(n)$ if $n
\in A$ and is $\eta^{\text{up}}_0(n)$ if $n \in \omega \backslash A$.  So
clearly $\alpha < \kappa_{\text{up}} \Rightarrow 
\neg(\eta^{\text{up}}_\alpha <_{J^{\text{bd}}_\omega} \eta')$
hence $\alpha < \kappa_{\text{up}} \Rightarrow
\neg(\eta^{\text{up}}_\alpha \le_{J^{\text{bd}}_\omega} \eta)$ hence by
clause $(\delta)$ of \ref{pt.0n}
for some $\beta < \kappa_{\text{dn}}$ we have
$\neg(\eta^{\text{dn}}_\beta <_{J^{\text{bd}}_\omega} \eta')$.  As
$\eta^{\text{dn}}_\beta <_{J^{\text{bd}}_\omega} \eta^{\text{up}}_0$,
necessarily $\neg(\eta^{\text{dn}}_\beta <_{J^{\text{bd}}_A} \eta')$
but $\gamma \in [\beta,\kappa_{\text{dn}}) \Rightarrow
\eta^{\text{dn}}_\beta \le_{J^{\text{bd}}_A} \eta^{\text{dn}}_\gamma$
hence $\gamma \in [\beta,\kappa_{\text{dn}}) \Rightarrow
\neg(\eta^{\text{dn}}_\gamma <_{J^{\text{bd}}_A} \eta') 
\Rightarrow \neg(\eta^{\text{dn}}_\gamma <_{J^{\text{bd}}_A} \eta)$,
as required.
\end{PROOF}

\noindent
We need the following from \cite[2.1]{Sh:885}:
\begin{claim}
\label{pt.885}  
Assume that $\kappa_1 \le \kappa_2$
are infinite regular cardinals, and there exists a
$(\kappa_1,\kappa_2)$-peculiar cut in ${}^\omega \omega$.

\Then \, for some $\sigma$-centered forcing notion $\bbQ$ of
cardinality $\kappa_1$ and a sequence $(I_\alpha:\alpha < \kappa_2)$
of open dense subsets of $\bbQ$, there is no directed $G \subseteq
\bbQ$ such that $(\forall \alpha < \kappa_2)(G \cap I_\alpha \ne
\emptyset)$.  Hence {\rm MA}$_{\kappa_2}$ fails.
\end{claim}

\begin{theorem}
\label{pt.1}  
Assume $\lambda = \text{\rm cf}(\lambda) =
\lambda^{<\lambda} > \aleph_2,\lambda > \kappa = \text{\rm cf}(\kappa) \ge
\aleph_1$ and $2^\lambda = \lambda^+$ and $(\forall \mu <
\lambda)(\mu^{\aleph_0} < \lambda)$.

For some forcing $\bbP^*$ of cardinality $\lambda^+$ not adding new
members to ${}^\lambda \bold V$ and $\bbP$-name $\name{\bbQ}^*$ 
of a c.c.c. forcing we have
$\Vdash_{\bbP^* * \name{\bbQ}^*} ``2^{\aleph_0}
= \lambda^+$ and ${\gp} = \lambda$ and {\rm MA}$_{< \lambda}$ and there
is a pair $(\bar f_1,\bar f^1)$ which is a peculiar
$(\kappa,\lambda)$-cut".
\end{theorem}

\begin{remark}
\label{pt.10} 
1) The proof of \ref{pt.1} is done in \S4 and broken into a 
series of Definitions and Claims, in particular we specify some of the
free choices in the general iteration theorem.

\noindent
2) In \ref{pt.3}(1), is cf$(\delta) > \aleph_0$ necessary?

\noindent
3) What if $\lambda = \aleph_2$?  The problem is \ref{p.7}(2).  To
eliminate this we may, instead quoting \ref{p.7}(2), start by
forcing $\bar \eta = \langle \eta_\alpha:\alpha <
\omega_1\rangle$ in $\bbP_{\bold k_0}$ and change some points.
\end{remark}

\noindent
Complementary to \ref{pt.1} is
\begin{observation}
\label{up.84} 
Assume $\lambda = \text{ cf}(\lambda)
> \aleph_1$ and $\mu = \text{ cf}(\mu) = \mu^{<\lambda} > \lambda$
\then \, for some c.c.c. forcing notion $\bbP$ of cardinality $\mu$
we have:

$\Vdash_{\bbP} ``2^{\aleph_0} = \mu,{\gp} = \lambda$ and for no
regular $\kappa < \lambda$ is there a peculiar $(\kappa,\lambda)$-cut
so ${\gt} = \lambda"$.
\end{observation}

\begin{PROOF}{\ref{up.84}}
We choose $\bar{\bbQ} = \langle \bbP_\alpha,\name{\bbQ}_\beta:\alpha \le
\mu,\beta < \mu\rangle$ such that:
\mn
\begin{enumerate}
\item[$(a)$]   $\bar{\bbQ}$ is an FS-iteration
\sn
\item[$(b)$]   $\name{\bbQ}_\beta$ is a
$\sigma$-centered forcing notion of cardinality $< \lambda$
\sn
\item[$(c)$]  if $\alpha < \mu,\name{\bbQ}$ is
a $\bbP_\alpha$-name of a $\sigma$-centered forcing notion of
cardinality $< \lambda$ \then \, for some $\beta \in [\alpha,\mu)$ we
have $\name{\bbQ}_\beta = \name{\bbQ}$
\sn
\item[$(d)$]   $\bbQ_0$ is adding $\lambda$ Cohens, $\langle
\name r_\varepsilon:\varepsilon <\lambda\rangle$.
\end{enumerate}
\mn
Clearly in $\bold V^{\bbP_\lambda}$ we have 
$2^{\aleph_0} = \lambda$, also every 
$\sigma$-centered forcing notion of cardinality $< \mu$, is from
$\bold V^{\bbP_\alpha}$ for some $\alpha < \mu$, so
as $\mu$ is regular we have
\mn
\begin{enumerate}
\item[$(*)$]  MA for $\sigma$-centered forcing notions of cardinality
$\le \lambda$ or just $< \mu$ dense sets
\end{enumerate}
\mn
Hence by \ref{pt.885} there is no peculiar
$(\kappa_1,\kappa_2)$-cut when $\aleph_1 \le \kappa_1 < \kappa_2 =
\lambda$ (even $\kappa_1 < \kappa_2 < \mu,\kappa_1 <
\lambda < \mu$). 
\end{PROOF}
\newpage

\section {Some specific forcing} 

\begin{definition}
\label{p.1}  
Let $\bar \eta =: \langle \eta_\alpha:\alpha < \alpha^*\rangle$ be a 
sequence of members of ${}^\omega \omega$ which is  
$<_{J^{\text{bd}}_\omega}$-increasing or just 
$\le_{J^{\text{bd}}_\omega}$-directed.  
We define the set ${\cF}_{\bar \eta}$ and
the forcing notion $\bbQ = \bbQ_{\bar \eta}$ and a generic real
$\name \nu$ for $\bbQ = \bbQ_{\bar \eta}$ as follows:
\mn
\begin{enumerate}
\item[$(a)$]   ${\cF}_{\bar \eta} = \{\nu \in {}^\omega(\omega +1)$:
if $\alpha < \ell g(\bar \eta)$ then $\eta_\alpha 
<_{J^{\text{bd}}_\omega}\nu\}$, here $\bar\eta$ is not\footnote{it is
enough that $\bar\eta$ is $\aleph_2$-directed by
$<_{J^{\text{bd}}_\omega}$; assuming $\bar\eta$ is just directed we
have to change clause $(e)(\beta)$ to $\eta_{\alpha(p)}
\le_{J^{\text{bd}}_\omega} \eta_{\alpha(q)}$, the situation}
necessarily $<_{J^{\text{bd}}_\omega}$-increasing
\sn
\item[$(b)$]   $\bbQ$ has the set of elements consisting of all
triples $p = (\rho,\alpha,g) =
(\rho^p,\alpha^p,g^p)$ (and $\alpha(p) = \alpha^p)$ such that
\begin{enumerate}
\item[$(\alpha)$]   $\rho \in {}^{\omega >}\omega$,
\sn
\item[$(\beta)$]  $\alpha < \ell g(\bar \eta)$, 
\sn
\item[$(\gamma)$]  $g \in {\cF}_{\bar\eta}$, and
\sn
\item[$(\delta)$]  if $n \in [\ell g(\rho),\omega)$ then
$\eta_\alpha(n) \le g(n)$; 
\end{enumerate}
\item[$(c)$]   $\le_{\bbQ}$ is defined by: $p \le_{\bbQ} q$ iff
(both are elements of $\bbQ$ and)
\begin{enumerate}
\item[$(\alpha)$]  $\rho^p \trianglelefteq \rho^q$,
\sn
\item[$(\beta)$]   $\alpha^p \le \alpha^q,\eta_{\alpha^p}
\le_{J^{\text{bd}}_\omega} \eta_{\alpha^q}$ so if $\bar\eta$ is
$<_{J^{\text{bd}}_\omega}$-increasing this means
\sn
\item[$(\gamma)$]  $g^q \le g^p$,
\sn
\item[$(\delta)$]   if $n \in [\ell g(\rho^q),\omega)$ then
$\eta_{\alpha(p)}(n) \le \eta_{\alpha(q)}(n)$,
\sn
\item[$(\varepsilon)$]  if $n \in [\ell g(\rho^p),\ell
g(\rho^q))$ then $\eta_{\alpha(p)}(n) \le \rho^q(n) \le g^p(n)$.
\end{enumerate}
\item[$(d)$]   For ${\cF} \subseteq {\cF}_{\bar \eta}$ 
which is downward directed (by $<_{J^{\text{bd}}_\omega}$)
we  define $\bbQ_{\bar \eta,{\cF}}$ as $\bbQ_{\bar \eta}
 \restriction \{p \in \bbQ_{\bar\eta}:g^p \in {\cF}\}$
\sn
\item[$(e)$]   $\name \nu = 
\name \nu_{\bbQ} = \name \nu_{\bbQ_{\bar \eta}} = \cup\{\rho^p:p \in
\name G_{\bbQ_{\bar \eta}}\}$.
\end{enumerate}
\end{definition}

\begin{claim}
\label{p.7}  
1) If $\bar \eta \in {}^\gamma({}^\omega
\omega)$ \then \, ${\cF}_{\bar \eta}$ is downward directed, 
in fact if $g_1,g_2 \in {\cF}_{\bar \eta}$ 
then $g = \text{\rm min}\{g_1,g_2\} \in {\cF}_{\bar \eta}$,
i.e., $g(n) = \text{\rm min}\{g_1(n),g_2(n)\}$ for $n < \omega$.  
Also ``$f \in {\cF}_{\bar \eta}"$ is absolute.

\noindent
[But possibly for every $\nu \in {}^\omega(\omega +1)$ we have: $\nu
\in {\cF}_{\bar\eta} \Leftrightarrow (\forall^* n)\nu(n) = \omega$].

\noindent
2) If $\bar \eta \in {}^\delta({}^\omega \omega)$ is 
$<_{J^{\text{bd}}_\omega}$-increasing and 
{\rm cf}$(\delta) > \aleph_1$ \then \, $\bbQ_{\bar \eta}$ is c.c.c.

\noindent
3) Moreover any set of $\aleph_1$ members of $\bbQ_{\bar \eta}$ is
included in the union of countably many directed subsets of 
$\bbQ_{\bar \eta}$.

\noindent
4) Assume $\langle \bbP_\varepsilon:\varepsilon \le \zeta\rangle$ is a
$\lessdot$-increasing sequence of c.c.c. forcing notions,
$\name{\bar \eta} = \langle
\name \eta_\alpha:\alpha < \delta\rangle$ is a 
$\bbP_0$-name of a $<_{J^{\text{bd}}_\omega}$-increasing sequence of
members of ${}^\omega \omega$ and {\rm cf}$(\delta) > \aleph_1$.  For
$\varepsilon \le \zeta$ let $\name{\bbQ}_\varepsilon$ be 
the $\bbP_\varepsilon$-name of the forcing
 notion $\bbQ_{\name{\bar \eta}}$ as defined in
$\bold V^{\bbP_\varepsilon}$.  \Then \, $\Vdash_{\bbP_\zeta} 
``\name {\bbQ}_\varepsilon$ is
$\subseteq$-increasing and $\le_{\text{ic}}$-increasing for $\varepsilon
\le \zeta$ and it is c.c.c. and {\rm cf}$(\zeta) > \aleph_0 \Rightarrow
\name{\bbQ}_\zeta = \cup\{\name{\bbQ}_\varepsilon:
\varepsilon < \zeta\}$ is c.c.c."

\noindent
5) Let $\bar \eta \in {}^\delta({}^\omega \omega)$ be as in part (2).
\mn
\begin{enumerate}
\item[$(a)$]   If ${\cF} \subseteq {\cF}_{\bar \eta}$ is
downward directed (by $\le_{J^{\text{bd}}_\omega}$) then 
$\bbQ_{\bar \eta,{\cF}}$ is  absolutely c.c.c.
\sn
\item[$(b)$]    If ${\cF}_1 \subseteq {\cF}_2 \subseteq 
{\cF}_{\bar \eta}$ are downward directed then $\bbQ_{\bar
\eta,{\cF}_1} \le_{\text{ic}} \bbQ_{\bar \eta,{\cF}_2}$.
\end{enumerate}
\mn
6) 
\mn
\begin{enumerate}
\item[$(a)$]   $\Vdash_{\bbQ_{\bar \eta}} ``\name \nu \in 
{}^\omega \omega$ and $\bold V[\name G] = 
\bold V[\name \nu]"$
\sn
\item[$(b)$]   $p \Vdash_{\bbQ_{\bar \eta}} 
``\rho^p \triangleleft \name \nu$ and $n \in 
[\ell g(\rho),\omega) \Rightarrow \eta_{\alpha(p)}(n) \le
\name \nu(n) \le g^p(n)"$
\sn
\item[$(c)$]   $\Vdash_{\bbQ_{\bar \eta}} ``p \in G$ iff $\rho
\triangleleft \name \nu \wedge (\forall n)(\ell
g(\rho) \le n < \omega \Rightarrow \eta_\alpha(n) \le \name \nu(n) 
\le g^p(n))"$
\sn
\item[$(d)$]   $\Vdash_{\bbQ_{\bar \eta}} 
``\name \nu \in {\cF}_{\bar \eta}$, i.e. 
$\name \nu(n) \in {\cF}^{\bold V[\bbQ_{\bar\eta}]}"$
\sn
\item[$(e)$]   $\Vdash_{\bbQ_{\bar \eta}}$ ``for every $f \in
({}^\omega \omega)^{\bold V}$ we have $f \in {\cF}_{\bar \eta}$ iff
$f \in {\cF}^{\bold V}_{\bar \eta}$ iff $\name \nu
\le_{J^{\text{bd}}_\omega} f"$.
\end{enumerate}
\end{claim}

\begin{PROOF}{\ref{p.7}}
1) Trivial.

\noindent
2) Assume $p_\varepsilon \in \bbQ_{\bar \eta}$ for $\varepsilon <
\omega_1$.  So $\{\alpha(p_\varepsilon):\varepsilon < \omega_1\}$ is
a set of $\le \aleph_1$ ordinals $< \delta$.  But cf$(\delta) >
\aleph_1$ hence there is
$\alpha(*) < \delta$ such that $\varepsilon < \omega_1 \Rightarrow
\alpha(p_\varepsilon) < \alpha(*)$.  For each $\varepsilon$ let
$n_\varepsilon = \text{ Min}\{n$: for every $k \in [n,\omega)$ we have
$\eta_{\alpha(p_\varepsilon)}(k) \le \eta_{\alpha(*)}(k) \le
g^{p_\varepsilon}(k)\}$.  It is well defined because
$\eta_{\alpha(p_\varepsilon)} <_{J^{\text{bd}}_\omega}
\eta_{\alpha(*)} <_{J^{\text{bd}}_\omega} g^{p_\varepsilon}$ recalling
$\alpha(p_\varepsilon) < \alpha(*)$ and $g^{p_\varepsilon} 
\in {\cF}_{\bar \eta}$.

So clearly for some 
$\bold x = (\rho^*,n^*,\eta^*,\nu^*)$ the following set is
uncountable

\begin{equation*}
\begin{array}{clcr}
{\cU} = {\cU}_{\bold x} = \{\varepsilon < \omega_1:
&\rho^{p_\varepsilon} = \rho^* \text{ and } n_\varepsilon
= n^* \text{ and }
\eta_{\alpha(p_\varepsilon)} \restriction n^*=\eta^* \\
  &\text{and } g^{p_\varepsilon} \restriction n^* = \nu^*\}.
\end{array}
\end{equation*}
\mn
Let

\begin{equation*}
\begin{array}{clcr}
\bbQ' = \bbQ'_{\bold x} =: 
\{p \in \bbQ_{\bar \eta}:&\ell g(\rho^p) \ge \ell g(\rho^*),
\rho^p \restriction \ell g(\rho^*) = \rho^*
\text{ and } \rho^p \restriction [\ell g(\rho^*),\ell g(\rho^p)) \subseteq
\eta_{\alpha(*)} \\
  &\text{ and } \alpha(p) < \alpha(*), \text{ and } \eta_{\alpha(p)}
  \restriction n^* = \eta^* \text{ and } g^p \restriction n^* = \nu^* \\
  &\text{ and } n \in [n^*,\omega) \Rightarrow \eta_{\alpha(p)}(n)
  \le \eta_{\alpha(*)}(n) \le g^p(n)\}.
\end{array}
\end{equation*}
\mn
Clearly
\mn
\begin{enumerate}
\item[$\circledast_1$]  $\{p_\varepsilon:\varepsilon \in {\cU}\}
 \subseteq \bbQ'$
\sn
\item[$\circledast_2$]  $\bbQ' \subseteq \bbQ_{\bar \eta}$ is
directed.
\end{enumerate}
\mn
So we are done.

\noindent
3) The proof of part (2) proves this.

\noindent
4),5) First we can check clause (b) of part (5) by the definitions of 
$\bbQ_{\bar \eta,{\cF}},\bbQ_{\bar \eta}$.
Second, concerning ``$\bbQ_{\bar \eta,{\cF}}$ is absolutely
c.c.c." (i.e. clause (a) of part (5)) note that if $\bbP$ is c.c.c.,
$G \subseteq \bbP$ is generic over $\bold V$ \then \, 
$\bbQ^{\bold V}_{\bar \eta,{\cF}} = 
\bbQ^{\bold V[G]}_{\bar \eta,{\cF}}$ and 
$\bbQ^{\bold V}_{\bar\eta,{\cF}} \le_{\text{ic}} 
\bbQ_{\bar \eta}^{\bold V} \le_{\text{ic}} 
\bbQ^{\bold V[G]}_{\bar\eta}$ by clause (b) and the last one is c.c.c. (as
$\bold V[G] \models ``\text{cf}(\ell g(\bar \eta)) > \aleph_1"$).  Hence 
$\bbQ^{\bold V}_{\bar \eta,{\cF}}$ is c.c.c. even in $\bold V[G]$
as required.  Turning to part (4), 
letting ${\cF}_\varepsilon = 
({\cF}_{\bar \eta})^{\bold V[\bbP_\varepsilon]}$, 
clearly $\Vdash_{\bbP_{\varepsilon_2}} 
``\name{\bbQ}_{\varepsilon_1} = \name{\bbQ}_{\bar \eta,
{\cF}_{\varepsilon_1}}"$ for $\varepsilon_1 <
 \varepsilon_2 < \zeta$.  Now about the c.c.c., as 
$\bbP_\varepsilon$ is c.c.c., it preserves ``cf$(\delta) > \aleph_1"$, so
 the proof of part (1) works.

\noindent
6) Easy, too.
\end{PROOF}

\begin{definition}
\label{p.14}  
Assume $\bar A = \langle
A_\alpha:\alpha < \alpha^*\rangle$ is a $\subseteq^*$-decreasing
sequence of members of $[\omega]^{\aleph_0}$.  We define the forcing
notion $\bbQ_{\bar A}$ and the generic real $\name w$ by:
\mn
\begin{enumerate}
\item[$(A)$]  $p \in \Bbb Q_{\bar A}$ \underline{iff}
\begin{enumerate}
\item[$(a)$]  $p = (w,n,A_\alpha) = (w_p,n_p,A_{\alpha(p)})$, 
\sn
\item[$(b)$]  $w \subseteq \omega \text{ is finite}$,
\sn
\item[$(c)$]  $\alpha < \alpha^* \text{ and } n < \omega$,
\end{enumerate}
\item[$(B)$]  $p \le_{\Bbb Q_{\bar A}} q$ \underline{iff}
\begin{enumerate}
\item[$(a)$]  $w_p \subseteq w_q \subseteq w_p \cup 
(A_{\alpha(p)} \backslash n_p)$
\sn
\item[$(b)$]  $n_p \le n_q$
\sn
\item[$(c)$]  $A_{\alpha(p)} \backslash n_p \supseteq A_{\alpha(q)}
  \backslash n_q$
\end{enumerate}
\item[$(C)$]  $\name w = \cup\{w_p:p \in \name G_{\bbQ_{\bar A}}\}$.
\end{enumerate}
\end{definition}

\begin{claim}
\label{p.15}  
Let $\bar A$ be as in Definition \ref{p.14}.

\noindent
1) $\bbQ_{\bar A}$ is a c.c.c. and even $\sigma$-centered forcing notion.

\noindent
2) $\Vdash_{\bbQ_{\bar A}} ``\name w \in [\omega]^{\aleph_0}$ 
is $\subseteq^* A_\alpha$ for each $\alpha <
\alpha^*"$ and $\bold V[\name G] = \bold V[\name w]$.

\noindent
3) Moreover, for every $p \in \bbQ_{\bar A}$ we have $\Vdash ``p \in
\name G$ iff $w_p \subseteq w \subseteq 
(A_{\alpha(p)} \backslash n_p) \cup w_p"$.
\end{claim}

\begin{PROOF}{\ref{p.15}}
  Easy.
\end{PROOF}

\begin{claim}
\label{p.21}  
Assume $\bar \eta \in {}^\delta({}^\omega
\omega)$ is $\le_{J^{\text{bd}}_\omega}$-increasing.

\noindent
1) If ${\cF} \subseteq {\cF}_{\bar \eta}$ is downward cofinal in
$({\cF}_{\bar \eta},<_{J^{\text{bd}}_\omega})$, i.e. $(\forall
 \nu \in {\cF}_{\bar\eta})(\exists \rho \in {\cF})(\rho
 <_{J^{\text{bd}}_\omega} \nu)$ and ${\cU} \subseteq \delta$ is
 unbounded then $\bbQ_{\bar \eta \restriction {\cU},{\cF}} =
 \{p \in \bbQ_{\bar \eta}:\alpha^p \in {\cU}$ and $g^p \in
 {\cF}\}$ is (not only $\subseteq \bbQ_{\bar\eta}$ but also is)
 a dense subset of $\bbQ_{\bar\eta}$.

\noindent
2) If {\rm cf}$(\delta) > \aleph_0$ and $\bbR$ is Cohen forcing then
$\Vdash_{\bbR} ``\bbQ^{\bold V}_{\bar\eta}$ is dense in 
$\bbQ^{\bold V[\name G]}_{\bar\eta}"$. 
\end{claim}

\begin{remark}
1) We can replace ``$\eta_\alpha
\le_{J^{\text{bd}}_\omega} \rho$" by ``$\rho$ belongs to the
$F_\sigma$-set $\bold B_\alpha$", where $\bold B_\alpha$ denotes a
Borel set from the ground model, i.e. its definition.

\noindent
2) Used in \ref{pg.10}.
\end{remark}

\begin{PROOF}{\ref{p.21}}
1) Check.

\noindent
2) See next claim. 
\end{PROOF}

\begin{claim}
\label{p.28} 
Let $\bar \eta = \langle
\eta_\gamma:\gamma < \delta\rangle$ is
$\le_{J^{\text{bd}}_\omega}$-increasing in ${}^\omega \omega$.

\noindent
1)  If $\bbP$ is a forcing notion of cardinality $< \text{\rm cf}
(\delta)$ \then \, $\Vdash_{\bbP} ``\bbQ^{\bold V}_{\bar
\eta}$ is dense in $\bbQ^{\bold V[\name G_{\bar \eta}]}"$.

\noindent
2) A sufficient condition for the conclusion of part (1) is:
\mn
\begin{enumerate}
\item[$\odot^{\text{cf}(\delta)}_{\Bbb P}$]  for every 
$X \in [\bbP]^{\text{cf}(\delta)}$ there is 
$Y \in [\bbP]^{<\text{cf}(\delta)}$
\newline
such that $(\forall p \in X)(\exists q \in Y)(p \le q)$.
\end{enumerate}
\mn
2A) We can weaken the condition to: if $X \in 
[\bbP]^{\text{\rm cf}(\delta)}$ then for some $q \in \bbP$, 
{\rm cf}$(\delta) \le |\{p \in X:p \le_{\bbP} q\}|$.

\noindent
3) If $\langle A_\alpha:\alpha < \delta^*\rangle$ is
$\subseteq^*$-decreasing sequence of infinite subsets of $\omega$ and
{\rm cf}$(\delta^*) \ne \text{\rm cf}(\delta)$ \then \,
$\odot^{\text{cf}(\delta)}_{\bbQ_{\bar A}}$ holds.
\end{claim}

\begin{PROOF}{\ref{p.28}}
1) By part (2).

\noindent
2) Let ${\cU} \subseteq \delta$ be unbounded of order type
cf$(\delta)$.  Assume $p \in \bbP$ and $\name \nu$
satisfies $p \Vdash_{\bbP} ``\name \nu \in 
{\cF}^{\bold V[\name G]}_{\bar\eta}"$.  So for every
$\gamma \in {\cU}$ we have $p \Vdash_{\bbP} ``\eta_\gamma
<_{J^{\text{bd}}_\omega} \name \nu \in {}^\omega
\omega"$, hence there is a pair $(p_\gamma,n_\gamma)$ such that:
\mn
\begin{enumerate}
\item[$(*)$]   $(a) \quad p \le_{\bbP} p_\gamma$
\sn
\item[${{}}$]   $(b) \quad n_\gamma < \omega$
\sn
\item[${{}}$]  $(c) \quad p_\varepsilon \Vdash_{\bbP} ``(\forall
n)(n_\gamma \le n < \omega \Rightarrow \eta_\varepsilon(n) <
\name \nu(n))$.
\end{enumerate}
\mn
We apply the assumption to the set 
$X = \{p_\varepsilon:\gamma \in {\cU}\}$
and get $Y \in [\bbP]^{< \text{cf}(\delta)}$ as there.  So for
every $\gamma \in {\cU}$ there is $q_\gamma$ such that $p_\gamma
\le_{\bbP} q_\gamma \in Y$.  As $|Y \times \omega| = |Y| + \aleph_0
< \text{ cf}(\delta) = |{\cU}|$ there is a pair $(q_*,n_*) \in Y
\times \omega$ such that ${\cU}' \subseteq \delta$ is unbounded
where ${\cU}' := \{\gamma \in {\cU}:q_\gamma = q_*$ and $n_\gamma
= n_*\}$.  Lastly, define $\nu_* \in {}^\omega(\omega +1)$ by
$\nu_*(n)$ is $0$ if $n < n_*$ is $\cup\{\eta_\alpha(n)+1:\alpha \in
{\cU}'\}$ when $n \ge n_*$.

Clearly
\mn
\begin{enumerate}
\item[$\circledast$]   $(a) \quad \nu_* \in {}^\omega(\omega +1)$
\sn
\item[${{}}$]   $(b) \quad \gamma \in {\cU}' \Rightarrow
\eta_\alpha \restriction [n_*,\omega) < \nu_* \restriction
[n_*,\omega)$
\sn
\item[${{}}$]   $(c) \quad$ if $\gamma < \delta$ then $\eta_\alpha
<_{J^{\text{bd}}_\omega} \nu_*$
\sn 
\item[${{}}$]   $(d) \quad \nu_* \in {\cF}^{\bold V}_{\bar\eta}$
\sn
\item[${{}}$]   $(e) \quad p \le q_*$
\sn
\item[${{}}$]   $(f) \quad q_* \Vdash_{\bbP} ``\nu_* \le
\name \nu"$.
\end{enumerate}
\mn
So we are done.

\noindent
2A) Similarly.

\noindent
3) If cf$(\delta^*) < \text{ cf}(\delta)$ let ${\cU} \subseteq
\delta^*$ be unbounded of order type cf$(\delta^*)$ and $\bbQ'_{\bar
A} = \{p\in \bbQ_{\bar A}:\alpha^p \in {\cU}\}$, it is dense in
$\bbQ_{\bar A}$ and has cardinality $\le \aleph_0 + 
\text{ cf}(\delta^*) < \text{ cf}(\delta)$, so we are done.

If cf$(\delta^*) > \text{ cf}(\delta)$ and $X \in 
[\bbP]^{\text{cf}(\delta)}$, let $\alpha(*) = \sup\{\alpha^p:p \in X\}$
and $Y = \{p \in \bbQ_{\bar A}:\alpha^p = \alpha(*)\}$.
\newline
The rest should be clear. 
\end{PROOF}
\newpage

\section {Proof of Theorem \ref{pt.1}}

\begin{choice}
\label{pt.3} 
1) $S \subseteq \{\delta < \lambda:
\text{cf}(\delta) > \aleph_0\}$ stationary.

\noindent
2) $\bar \eta$ is as in \ref{pt.7} below, so possibly a preliminary
forcing of cardinality $\aleph_2$ we have such $\bar \eta$ with
cf$(\ell g(\bar \eta)) > \aleph_1$.
\end{choice}

\begin{defclaim}
\label{pt.7}
1) Assume $\kappa = \text{ cf}(\kappa) \in 
[\aleph_2,\lambda)$ and $\bar \eta = \langle 
\eta_\alpha:\alpha < \kappa\rangle$ is an
$<_{J^{\text{bd}}_\omega}$-increasing sequence in ${}^\omega \omega$
and $\delta \in \lambda \backslash \omega_1$ a limit ordinal and $\gamma
\le \lambda$.  \Then \, the following $\bold s = \bold
s_{\delta,\gamma}$ is a winning strategy
of COM in the game $\Game_{\delta,< \gamma}$: COM just preserves:
\mn
\begin{enumerate}
\item[$\otimes$]   $(a) \quad$ if for every $\zeta < \varepsilon$ we
have $(\alpha) + (\beta)$ \then \, we have $(*)$ where
\begin{enumerate}
\item[${{}}$]  $(\alpha) \quad \bbP_{\bold q_\zeta} = \bbP_{\bold p_\zeta} * 
\name{\bbQ}_{\bar \eta}$ where $\name{\bbQ}_{\bar \eta}$ 
is from \ref{p.1} and in $\bold V^{\bbP[\bold p_\zeta]}$, i.e. is a

\hskip25pt $\bbP_{\bold p_\zeta}$-name
\sn
\item[${{}}$]  $(\beta) \quad \bbP_{\bold p_\zeta} * 
\name{\bbQ}_{\bar \eta} \lessdot \bbP_{\bold p_\varepsilon} *
\name{\bbQ}_{\bar \eta}$ 
\sn
\item[${{}}$]  $(*) \quad \bbP_{\bold q_\varepsilon} 
= \bbP_{\bold p_\zeta} * \name{\bbQ}_{\bar \eta}$,  
so we have to interpret $\bbP_{\bold q_\varepsilon}$ such that its set
of

\hskip35pt  elements is $\subseteq {\cH}_{< \aleph_1}
(u^{\bold q_\varepsilon})$ which is easy, i.e. it is
$\bbP_{\bold p_\varepsilon} \cup \{(p,\name r)$:

\hskip35pt $p \in \bbP_{\bold p_\varepsilon}$ and $\name r$ is a canonical
$\bbP_{\bold p_\varepsilon}$-name of a 
member of $\name {\bbQ}_{\bar \eta}$

\hskip35pt  (i.e. use $\aleph_0$ maximal antichains, etc.)$\}$
\end{enumerate}
\item[${{}}$]  $(b) \quad$ if in (a) clause $(\alpha)$ holds but
$(\beta)$ fail \then
\begin{enumerate}
\item[${{}}$]  $(\alpha) \quad$ the set of elements of 
$\bbP_{\bold q_\varepsilon}$ is $\bbP_{\bold p_\varepsilon} \cup
\{(p,\name r)$: for some $\zeta < \varepsilon$ and

\hskip25pt $(p',\name r) \in \bbP_{\bold q_\zeta}$ we have
$\bbP_{\bold p_\varepsilon} \models ``p' \le p"$ 
\sn
\item[${{}}$]  $(\beta) \quad$ the order is defined naturally
\end{enumerate}
\item[${{}}$]   $(c) \quad$ if in (a), clause $(\alpha)$ fail, let
$\zeta$ be minimal such that it fails, and then
\begin{enumerate}
\item[${{}}$]   $(\alpha) \quad$ the set of 
elements of $\bbP_{\bold q_\varepsilon}$ is
$\bbP_{\bold p_\varepsilon} \cup \{(p,\name r)$:
for some $\xi < \zeta$ and $p'$

\hskip25pt  we have 
$(p',\name r) \in \bbP_{\bold q_\zeta}$ and
$\bbP_{\bold p_\varepsilon} \models ``p' \le p"\}$
\sn
\item[${{}}$]  $(\beta) \quad$ the order is natural.
\end{enumerate}
\end{enumerate}
\end{defclaim}

\begin{remark}
In \ref{pt.7} we can combine clauses (b) and (c). 
\end{remark}

\begin{PROOF}{\ref{pt.7}}
By \ref{p.7} this is easy, see in particular
\ref{p.7}(4). 
\end{PROOF}

\noindent
Technically it is more 
convenient to use the (essentially equivalent) variant.
\begin{defclaim}
\label{pg.10}  
1) We replace $\bbP_{\bold q_\zeta} = \bbP_{\bold p_\zeta} * 
\bbQ_{\bar\eta}$ by $\bbP_{\bold q_\zeta} = \bbP_{\bold p_\zeta} * 
\bbQ_{\bar\eta,{\cF}_\zeta}$ where

\begin{equation*}
\begin{array}{clcr}
{\cF}_\zeta = \{\nu:&\text{ for some } \varepsilon \le
\zeta,\nu \in {\cF}^{\bold V^[\bbP[\bold p_\varepsilon]]}_{\bar\eta}
\text{ but} \\
  &\text{ for no } \xi < \varepsilon \text{ and } \nu_1 \in
\cF^{\bold V[\bbP[\bold p_\xi]]}_{\bar\eta} \text{ do we have} \\
  &\nu_1 \le_{J^{\text{bd}}_\omega} \nu\}.
\end{array}
\end{equation*}
\mn
2) No change by \ref{p.21}(1).
\end{defclaim}

\begin{remark}
\label{pt.14}  
In \ref{pt.7} we can use $\name{\bar\eta} =
 \langle \name \eta_\alpha:\alpha < \kappa\rangle$
say a $\bbP_{\bold k_0}$-name, but \then \, for the game 
$\Game_{\delta,f(\delta)}$ we better assume $\delta \in E_{\bold k_0}$
 and $\name{\bar\eta}$ is a $\bbP[\bold p^{\bold k}_\delta]$-name.
\end{remark}

\begin{defclaim}
\label{pt.21}  
1) Let $\bold k_* \in K^2_\lambda$ and $\name \nu_\alpha \,(\alpha <
\lambda)$ be chosen as follows:
\mn
\begin{enumerate}
\item[$(a)$]   $E_{\bold k_*} = \lambda$ and
$u[\bold p^{\bold k_*}_\alpha] = \omega_1 +
\alpha$ hence $u[\bold k_*] = \lambda$
\sn
\item[$(b)$]  $\bbP^{\bold k_*}_\alpha$ is $\lessdot$-increasing
continuous
\sn
\item[$(c)$]  $\bbP^{\bold k_*}_{\alpha +1} = 
\bbP^{\bold k_*}_\alpha * \name{\bbQ}_{\bar \eta}$ and
$\name \nu_\delta$ is the generic (for this copy) of
$\name {\bbQ}_{\bar \eta}$ where $\bar\eta$ is from \ref{pt.7}
\sn 
\item[$(d)$]  $S_{\bold k_*} = S$ (a stationary subset of
$\lambda$), $\delta \in S \Rightarrow \text{ cf}(\delta) > \aleph_0$
\sn
\item[$(e)$]   for each $\delta \in S_{\bold k_*},\bold s^{\bold
k_*}_\delta = \bold s_{\delta,\lambda}$ is from \ref{pt.7} or better
\ref{pg.10} 
\sn
\item[$(f)$]  $\bold g^{\bold k_*}_\delta$ is $\langle(\bold
p^{\bold k_*}_\delta,\bold p^{\bold k_*}_{\delta +1})\rangle$,
mv$(\bold g^{\bold k_*}_\delta) = 0$, only one move was done.
\end{enumerate}
\mn
2) If $\bold k_* \le_{K_2} \bold k$ then $\Vdash_{\bbP_{\bold k}}$
``the pair $(\langle \name \nu_\alpha:\alpha < \lambda\rangle,
\langle \eta_i:i < \kappa\rangle)$ is a
$(\lambda,\kappa)$-peculiar cut". 
\end{defclaim}

\begin{PROOF}{\ref{pt.21}}
Clear (by \ref{pt.7}).  
\end{PROOF}

\begin{definition}
\label{pt.28}  
Let $\bbP^*$ be the following forcing notion:
\mn
\begin{enumerate}
\item[$(A)$]   the members are $\bold k$ such that
\begin{enumerate}
\item[$(a)$]    $\bold k_* \le_{K_2} \bold k \in K^2_\lambda$
\sn  
\item[$(b)$]    $u[\bold k] = \cup\{u[\bold p^{\bold k}_\alpha]:
\alpha \in E_{\bold k}\}$ is an ordinal $<
\lambda^+$ (but of course $\ge \lambda$)
\sn  
\item[$(c)$]   $S_{\bold k} = S_{\bold k_*}$ and $\bold s^{\bold
k}_\delta = \bold s^{\bold k_*}_\delta$ for $\delta \in S_{\bold k}$
\end{enumerate}
\item[$(B)$]   the order: $\le_{K^2_\lambda}$.
\end{enumerate}
\end{definition}

\begin{definition}
\label{pt.35}  
We define the $\bbP^*$-name $\name{\bbQ}^*$ as

\[
\cup\{\bbP^{\bold k}_\lambda:\bold k \in \name{\bbG}_{\bbP^*}\} 
= \cup\{\bbP_{\bold p}[\bold p^{\bold k}_\alpha]:
\alpha \in E_{\bold k} \text{ and } \bold k \in
\name G_{\bbP^*}\}.
\]
\end{definition}

\begin{claim}
\label{pt.42}  
1) $\bbP^*$ has cardinality $\lambda^+$.

\noindent
2) $\bbP^*$ is strategically $(\lambda +1)$-complete hence add no
new member to ${}^\lambda \bold V$.

\noindent
3) $\Vdash_{\bbP^*} ``\name{\bbQ}^*$ is
c.c.c. of cardinality $\le \lambda^+$".

\noindent
4) $\bbP^* * \name{\bbQ}^*$ is a forcing notion
of cardinality $\lambda^+$ neither collapsing any cardinal nor
changing cofinalities.

\noindent
5) If $\bold k \in \bbP^*$ \then \, $\bold k \Vdash_{\bbP^*}
``\bbP_{\bold k} \lessdot \name{\bbQ}^*"$ hence
$\Vdash_{\bbP^*} ``\bbP_{\bold k_*} \lessdot 
\name{\bbQ}^*"$. 
\end{claim}

\begin{PROOF}{\ref{pt.42}}
1) Trivial.

\noindent
2) By claim \ref{it.49}.

\noindent
3) $\name G_{\bbP^*}$ is $(< \lambda^+)$-directed.

\noindent
4),5)  Should be clear. 
\end{PROOF}

\begin{claim}
\label{pt.45}  
If $\bold k \in \bbP^*$ and $G \subseteq \bbP_{\bold k}$ 
is generic over $\bold V$ \then \,
\mn
\begin{enumerate}
\item[$(a)$]   $\langle \name \nu_\alpha[G \cap
\bbP_{\bold k_*}]:\alpha < \lambda\rangle$ is
$<_{J^{\text{bd}}_\omega}$-decreasing 
and $i < \kappa \Rightarrow \eta_i
<_{J^{\text{bd}}_\omega} \name \nu_\alpha[G \cap
\bbP_{\bold k_*}]$, (this concerns $\bbP_{\bold k_*}$ only)
\sn
\item[$(b)$]    if $\rho \in ({}^\omega \omega)^{\bold V[G]}$ and $i
< \kappa \Rightarrow \eta_i <_{J^{\text{bd}}_\omega} \rho$ \then \,
for every $\alpha < \lambda$ large enough we have
$\name \nu_\alpha[G] <_{J^{\text{bd}}_\omega} \rho$
\sn
\item[$(c)$]    if $\rho \in ({}^\omega \omega)^{\bold V[G]}$ and $i
< \kappa \Rightarrow \eta_i \nleq_{J^{\text{bd}}_\omega} \rho$ \then \,
for every $\alpha < \lambda$ large enough we have
$\name \nu_\alpha[G] \nleq_{J^{\text{bd}}_\omega} \rho$.
\end{enumerate}
\end{claim}

\begin{PROOF}{\ref{pt.45}}
Should be clear. 
\end{PROOF}

\begin{claim}
\label{pt.49}  
1) If $\bold k \in \bbP^*$ and
$\name{\bbQ}$ is a $\bbP_{\bold k}$-name of a
c.c.c. forcing of cardinality $< \lambda$ and 
$\alpha \in E_{\bold k}$ and $\name{\bbQ}$ is a 
$\bbP[\bold p^{\bold k}_\alpha]$-name
\then \, for some $\bold k_1$ we have:
\mn
\begin{enumerate}
\item[$(a)$]   $\bold k \le_{K_1} \bold k_1 \in \bbP^*$
\sn
\item[$(b)$]   $\Vdash_{\bbP_{\bold k_1}}$ ``there is a subset of 
$\name{\bbQ}$ generic over $\bold V[G_{\bbP_{\bold k_1}} 
\cap \bbP[\bold p^{\bold k}_\alpha]]"$.
\end{enumerate}
\mn
2) In (1) if $\Vdash_{\bbP[\bold p^{\bold k}_\alpha] * 
\name {\bbQ}}$ ``there is $\rho \in {}^\omega 2$ not in $\bold
V[\name G_{\bbP_{\bold k}}]"$ \then \, $\Vdash_{\bbP_{\bold k_1}}$ 
``there is $\rho \in {}^\omega 2$ not in $\bold
V[\name G_{\bbP_{\bold k}}]"$.
\end{claim}

\begin{PROOF}{\ref{pt.49}}
1) By \ref{it.56}.

\noindent
2) By part (1) and clause $(\eta)$ of \ref{it.56}.   
\end{PROOF}

\begin{PROOF}{\ref{pt.1}} 
\underline{Proof of Theorem \ref{pt.1}}  
We force by $\bbP^* * \name {\bbQ}^*$ where $\bbP^*$ is defined in
\ref{pt.28} and the $\bbP^*$-name $\name{\bbQ}^*$ 
is defined in \ref{pt.35}.  By Claim \ref{pt.42}(4) we know
that no cardinal is collapsed and no cofinality is changed.  We know
that $\Vdash_{\bbP^* * \name{\bbQ}^*}
``2^{\aleph_0} \le \lambda^+"$ because $|\bbP^*| = \lambda^+$ and
$\Vdash_{\bbP^*} ``\name {\bbQ}^*$ has
cardinality $\le \lambda^+"$, so $\bbP^* * \name{\bbQ}^*$ 
has cardinality $\lambda^+$, see \ref{pt.42}(3),(4).  

Also $\Vdash_{\bbP^* * \name{\bbQ}}
``2^{\aleph_0} \ge \lambda^+"$ as by \ref{pt.42}(2) it suffices to
prove: for every $\bold k_1 \in \bbP^*$ there is $\bold k_2 \in 
\bbP^*$ such that $\bold k_1 \le_{K_2} \bold k_2$ and forcing by 
$\bbP_{\bold k_2}/\bbP_{\bold k_1}$ adds a real, which holds by
\ref{pt.49}(2). 

Lastly, we have to prove that $(\langle \eta_i:i <
\kappa\rangle,\langle \name \nu_\alpha:\alpha <
\lambda\rangle)$ is a peculiar cut.  In Definition \ref{pt.0n} clauses
$(\alpha),(\beta),(\gamma)$ holds by the choice of $\bold k_*$.  As
for clauses $(\delta),(\varepsilon)$
to check this it suffices to prove that for every $f \in {}^\omega
\omega$ they hold, so it is suffice to check it in any 
sub-universe to which $(\bar\eta,\bar \nu),f$ belong.  
Hence by \ref{pt.42}(1) it
suffices to check it in $\bold V^{\bbP_{\bold k}}$ for any $\bold k
\in \bbP^*$.  But this holds by \ref{pt.21}(2).
\end{PROOF}
\newpage

\section {Quite general applications} 

\begin{theorem}
\label{bt.7}  
Assume $\lambda = \text{\rm cf}(\lambda) =
\lambda^{<\lambda} > \aleph_2$ and 
$2^\lambda = \lambda^+$ and $(\forall \mu < \lambda)(\mu^{\aleph_0} 
< \lambda)$.  \Then \, for 
some forcing $\bbP^*$ of cardinality $\lambda^+$ not adding new
members to ${}^\lambda \bold V$ and $\bbP^*$-name $\name{\bbQ}^*$ of 
a c.c.c. forcing it is forced, i.e. 
$\Vdash_{\bbP^* * \name{\bbQ}^*}$ that $2^{\aleph_0} = \lambda^+$ and 
\mn
\begin{enumerate}
\item[$(a)$]  ${\gp} = \lambda$ and {\rm MA}$_{< \lambda}$
\sn
\item[$(b)$]   for every regular $\kappa \in (\aleph_1,\lambda)$
there is a $(\kappa,\lambda)$-peculiar cut $(\langle \eta^\kappa_i:i <
\kappa\rangle,\langle \nu^\kappa_\alpha:\alpha < \lambda\rangle)$
hence $\gp = \gt = lambda$
\sn
\item[$(c)$]   if $\bbQ$ is a (definition of $a$) Suslin c.c.c. forcing
notion defined by $\bar\varphi$ possibly with a real parameter from
$\bold V$, \then \, we can find a sequence $\langle 
\nu_{\bbQ,\name \eta,\alpha}:\alpha < \lambda\rangle$ which is positive for
$(\bbQ,\name \eta)$, see \cite{Sh:630}, e.g. {\rm non(null)} $=\lambda$
\sn
\item[$(d)$]    in particular ${\gb} = {\gd} = \lambda$.
\end{enumerate}  
\end{theorem}

\begin{remark}
\label{bt.9}
0) In clause (c) we can let $\bbQ$ be a c.c.c nep  forcing 
(see \cite{Sh:630}), with ${\gB},{\gC}$ of cardinality $\le \lambda$ 
and $\name \eta$ is a $\bbQ$-name of a real (i.e. member of ${}^\omega
2$).

\noindent
1) Concerning \ref{bt.7} as remarked earlier in
\ref{it.30}(1), if
we like to deal with Suslin forcing defined with a real parameter from
$\bold V^{P^* * \name{\bbQ}^+}$ and similarly for $\gB,\gC$
 we in a sense have to change/create new strategies.  We could start with
$\langle S_\alpha:\alpha < \lambda^+\rangle$ such that $S_\alpha
\subseteq \lambda,\alpha < \beta \Rightarrow |S_\alpha \backslash
S_\beta| < \lambda$ and $S_{\alpha +1} \backslash S_\alpha$ is a
stationary subset of $\lambda$.  But we can code this in the
strategies, do nothing till you know the definition of the forcing.

\noindent
2) We may like to strengthen \ref{bt.7} by demanding
\mn
\begin{enumerate}
\item[$(c)$]   for some $\bbQ$ as in clause (c) of \ref{bt.7},
MA$_{\bbQ}$ holds or even for a dense set of $\bold k_1 \in
\bbP^*$, see below, there is $\bold k_2 \in \bbP^*$ such that
$\bold k_1 \le_{K_2} \bold k_2$ and $\bbP_{\bold k_2}/\bbP_{\bold
k_1}$ is $\bbQ^{\bold V[\bbP_{\bold k_1}]}$.
\end{enumerate}
\mn
For this we have to restrict the family of $\bbQ$'s in clause (c)
such that those two families are orthogonal, i.e. commute.  Note,
however, that for Suslin c.c.c forcing this is rare, see
\cite{Sh:630}.

\noindent
3) This solves the second Bartoszynski
test problem, i.e. (B) of Problem \ref{h.7}.

\noindent
4) So $(\bar\varphi,\bbQ,\nu,\name \eta)$ in clause (c) of \ref{bt.7}
   satisfies
\mn
\begin{enumerate}
\item[$(a)$]  $\nu \in {}^\omega 2$
\sn
\item[$(b)$]  $\bar\varphi = (\varphi_0,\varphi_1,\varphi_2),
\Sigma_1$ formulas with the real parameter $\nu$ 
\sn
\item[$(c)$]  $\bbQ$ is the forcing notion defined by: 
\begin{enumerate}
\item[$\bullet$]  set of elements $\{\rho \in {}^\omega 2:\varphi_0[\rho]\}$
\sn
\item[$\bullet$]  quasi order $\le_{\bbQ} = \{(\rho_1,\rho_2):\rho_1,\rho_2 \in
{}^\omega 2,\varphi_1(\rho_1,\rho_2)\}$
\sn
\item[$\bullet$]  incompatibility in $\bbQ$ is defined by $\varphi_3$
\end{enumerate}
\item[$(d)$]  $\name\eta$ is a $\bbQ$-name of a real, i.e. $\langle
p_{n,k}:k \le \omega\rangle$ a (absolute) maximal antichain of
$\bbQ,\bold t_k = \langle \bold t_{n,k}:k < \omega\rangle,\bold
t_{k,n}$ a truth value.
\end{enumerate}
\end{remark}

\begin{proof}  The proof is like the proof of \ref{pt.1} so
essentially broken to a series of definitions and Claims.
\end{proof}

\begin{claim}
\label{bt.14}
\underline{Claim/Choice}:

Without loss of generality there
is a sequence $\langle S_\alpha:\alpha < \lambda^+\rangle$ such that:
\mn
\begin{enumerate}
\item[$(a)$]   $S_\alpha \subseteq S^\lambda_{\aleph_0}$ is stationary
\sn
\item[$(b)$]   if $\alpha < \beta$ then $S_\alpha \backslash
S_\beta$ is bounded (in $\lambda$)
\sn
\item[$(c)$]   $\diamondsuit_{S_{\alpha +1}\backslash S_\alpha}$
and $\diamondsuit_{S^\lambda_{\aleph_0} \backslash \cup\{S_\alpha:\alpha <
\lambda^+\}}$. 
\end{enumerate}
\end{claim}

\begin{proof}  E.g. by a preliminary forcing.
\end{proof}

\begin{definition}
\label{bt.28} 
Let $\Bbb P^*$ be the following forcing notion:
\mn
\begin{enumerate}
\item[$(A)$]   The members are $\bold k$ such that
\begin{enumerate}
\item[$(a)$]   $\bold k \in K^2_\lambda$
\sn  
\item[$(b)$]   $u[\bold k] = \cup\{u[\bold p^{\bold k}_\alpha]:
\alpha \in E_{\bold k}\}$ is an ordinal $< \lambda^+$ (but 
of course $\ge \lambda$)
\sn  
\item[$(c)$]    $S_{\bold k} \in \{S_\alpha:\alpha < \lambda^+\}$.
\end{enumerate}
\item[$(B)$]   The order: $\le_{K^2_\lambda}$.
\end{enumerate}
\end{definition}

\begin{definition}
\label{bt.35}  We define the $\bbP^*$-name $\name{\bbQ}^*$ as

\[
\cup\{\bbP^{\bold k}_\lambda:\bold k \in \name{\bbG}_{\bbP^*}\} 
= \cup\{\bbP[\bold p^{\bold k}_\alpha]:\alpha \in E_{\bold k} 
\text{ and } \bold k \in \name G_{\bbP^*}\}.
\]
\end{definition}

\begin{claim}
\label{bt.42}  As in \ref{pt.42}:

\noindent
1) $\bbP^*$ has cardinality $\lambda^+$.

\noindent
2) $\bbP^*$ is strategically $(\lambda +1)$-complete hence add no
 new member to ${}^\lambda \bold V$.

\noindent
3) $\Vdash_{\bbP^*} ``\name{\bbQ}^*$ is
c.c.c. of cardinality $\le \lambda^+$".

\noindent
4) $\bbP^* * \name{\bbQ}^*$ is a forcing notion
of cardinality $\lambda^+$ neither collapsing any cardinal nor
changing cofinalities.

\noindent
5) If $\bold k \in \bbP^*$ \then \, $\bold k \Vdash_{\bbP^*}
``\bbP_{\bold k} \lessdot \name{\bbQ}^*"$ hence
$\Vdash_{\bbP^*} ``\bbP_{\bold k_*} \lessdot 
\name{\bbQ}^*"$. 
\end{claim}

\begin{PROOF}{\ref{pt.42}}
1) Trivial.

\noindent
2) By claim \ref{it.49}.

\noindent
3) $\name G_{\bbP^*}$ is $(< \lambda^+)$-directed.

\noindent
4),5)  Should be clear.  
\end{PROOF}

\begin{claim}
\label{bt.45} 
Assume
\mn
\begin{enumerate}
\item[$(A)$]  $(a) \quad \bold k \in \bbP^*$
\sn
\item[${{}}$]  $(b) \quad S_{\bold k} = S_\alpha,\alpha < \lambda^+$
\sn
\item[${{}}$]  $(c) \quad \name \nu$ is a 
$\bbP^{\bold k}_\varepsilon$-name of a member of 
${}^\omega 2,\varepsilon < \kappa$
\sn
\item[${{}}$]  $(d) \quad \name{\bbQ}$ is a 
$\bbP_{\bold k_1}$-name of a c.c.c. Suslin forcing and
$\name \eta$ a $\name{\bbQ}$-name both 

\hskip25pt definable from $\name \nu$.
\end{enumerate}
\mn
\Then \, there is $\bold k_2$ such that
\mn
\begin{enumerate}
\item[$(B)$]   $(a) \quad \bold k_1 \le \bold k_2$
\sn
\item[${{}}$]   $(b) \quad S_{\bold k_2} = S_{\alpha +1}$
\sn
\item[${{}}$]  $(c) \quad$ if $\varepsilon \in S_{\alpha +1} \backslash
S_\alpha$ \then \, $\bbP^{\bold k_2}_{\varepsilon +1} = 
\bbP^{\bold k_2}_\varepsilon * \name{\bbQ}$ and
$\name \eta_\varepsilon$ is the copy of $\name\eta$
\sn
\item[${{}}$]  $(d) \quad$ if $\varepsilon \in S_{\alpha +1} \backslash
  S_\varepsilon$ \then \, the strategy 
{\bf st}$_\varepsilon$ is as in \ref{pt.7}, using 
$\name{\bbQ}$ instead of 

\hskip25pt $\name{\bbQ}_{\bar \eta}$.
\end{enumerate}
\end{claim}

\begin{PROOF}{\ref{pt.45}}
Straight.  
\end{PROOF}

\begin{claim}
\label{bt.49}
Like \ref{pt.49}:

\noindent
1) If $\bold k \in \bbP^*$ and $\name{\bbQ}$ is a $\bbP_{\bold k}$-name of a
c.c.c. forcing of cardinality $< \lambda$ and 
$\alpha \in E_{\bold k}$ and $\name{\bbQ}$ is a 
$\bbP[\bold p^{\bold k}_\alpha]$-name
\then \, for some $\bold k'$ we have:
\mn
\begin{enumerate}
\item[$(a)$]   $\bold k \le_{K_2} \bold k_1 \in \bbP^*$
\sn
\item[$(b)$]   $\Vdash_{\bbP_{\bold k_1}}$ ``there is a subset of 
$\name{\bbQ}$ generic over $\bold V[G_{\bbP_{\bold k_1}} \cap 
\bbP[\bold p^{\bold k}_\alpha]]$.
\end{enumerate}
\mn
2) In (1) if $\Vdash_{\bbP_{\bold k} * \name{\bbQ}}$ ``there is 
$\rho \in {}^\omega 2$ not in $\bold V[\name G_{\bbP_{\bold k}}]"$ 
then $\Vdash_{\bbP_{\bold k_1}}$ ``there is $\rho \in {}^\omega 2$ 
not in $\bold V[\name G_{\bbP_{\bold k}}]"$.
\end{claim}

\begin{PROOF}{\ref{pt.49}}
1) By \ref{it.56}.

\noindent
2) By part (1) and clause $(\eta)$ of \ref{it.56}.   
\end{PROOF}

\begin{PROOF}{\ref{bt.7}}
\underline{Proof of Theorem \ref{bt.7}}  

We force by $\bbP^* * \name{\bbQ}^*$ where $\bbP^*$ is defined in
\ref{bt.28} and the $\bbP^*$-name $\name{\bbQ}$ is defined in 
\ref{bt.35}.  By Claim \ref{bt.42}(4) we know
that no cardinal is collapsed and no cofinality is changed.  We know
that $\Vdash_{\bbP^* * \name{\bbQ}^*}
``2^{\aleph_0} \le \lambda^+"$ because $|\bbP^*| = \lambda^+$ and
$\Vdash_{\bbP^*} ``\name{\bbQ}^*$ has
cardinality $\le \lambda^+"$, so $\bbP^* * \name{\bbQ}^*$ has
cardinality $\lambda^+$, see \ref{bt.42}(3),(4).  

Also $\Vdash_{\bbP^* * \name{\bbQ}}
``2^{\aleph_0} \ge \lambda^+"$ as by \ref{pt.42}(2) it suffices to
prove: for every $\bold k_1 \in \bbP^*$ there is $\bold k_2 \in 
\bbP^*$ such that $\bold k_1 \le_{K_2} \bold k_2$ and forcing by 
$\bbP_{\bold k_2}/\bbP_{\bold k_1}$ add a real, which holds by
\ref{bt.49}(2).  Similarly 
$\Vdash_{\bbP * \name{\bbQ}^*} ``\text{MA}_{< \lambda}"$ even for 
$\lambda$ dense subsets by \ref{bt.49}(1) we have proved clause (a) of
\ref{bt.7}.

Clause (b) of \ref{bt.7} is proved as in the proof of \ref{pt.1},
$\bold k_*$ is above $\bold k_0$.

As for clause (c) 
we are given $\bold k_0$ and $\bbQ,\name \nu,\name \eta$ such that
$\name \nu$ is a $(\bbP^*_-,Q^*)$-name of a real and $\name{\bbQ}$ is
a Suslin c.c.c. forcing definable (say by $\bar\varphi_0$) from the
real $\name \nu$ and $\name\eta$ a $(\bbP^* * \name{\bbQ}^*)$-name of
$\name{\bbQ}$-name for $\name{\bbQ}$ of a real defined by $\aleph_0$
maximal antichain of $\name{\bbQ}$, absolutely of course.

As $\Vdash_{\bbP^*} ``\name{\bbQ}^*$ satisfies the c.c.c.", for some
$\bold k_1 \in \bbP^*$ above $\bold k_0$ and $\bbP_{\bold k_1}$-name
$\name \nu'$ of a member of ${}^{\lambda \ge} 2$ and $\name \eta'$ is a
$\bbP_{\bold k_1}$-name in $\name{\bbQ}_{\bar\varphi,\nu'}$ we have
$\bold k_1 \Vdash_{\bbP^*} ``\name\nu = \name \nu' \wedge \name \eta =
\name \eta'"$.

As $\bbP_{\bold k_1}$ satisfies the c.c.c. for some $\varepsilon <
\lambda,(\bold k_1,\varepsilon,\name \nu',\name{\bbQ}_{\name
\nu'},\name \eta')$ satisfies the assumptions on $(\bold
k,\varepsilon,\name \nu',\name eta')$ is as in \ref{bt.45} so there is
$\bold k_2$ and $\langle \name \eta_\alpha:\alpha \in S_{\alpha +1}
\backslash S_\alpha\rangle$ as there.  So $\bold k_0 \le \bold k_1 \le
\bold k_2$ and
\mn
\begin{enumerate}
\item[$(*)$]  if $\bold k_2 \le \bold k_3$ \then \, for a club of
$\zeta < \lambda,\name \nu'$ is a $\bbP^{\bold k_3}_\zeta$-name and
$\name \eta_\zeta$ is $(\name{\bbQ}_{\bar\varphi,\\bar
\nu'},\name\eta)$-generic over $\bold V^{\bbP_\zeta[\bold k_3]}$.
\end{enumerate}
\mn
This is clearly enough, so clause (a) of \ref{bt.7} holds.  For clause
(d) of \ref{bt.7}, first
Random real forcing is a Suslin c.c.c. forcing so non(null) $\le \lambda$
follows from clause (c) and non(nul) $\ge \lambda$ follows from clause
(a).  

Lastly, ${\gb} \ge \lambda$ by MA$_{< \lambda}$ and we know
${\gd} \ge {\gb}$.  As dominating real forcing = Hechler forcing
is a c.c.c. Suslin forcing so by clause (c) we have ${\gd} \le
\lambda$, together ${\gd} = {\gb} = \lambda$, i.e. clause (d)
holds. 
\end{PROOF}
\newpage


\end{document}